\documentclass[preprint,12pt]{elsarticle}

\usepackage{fullpage}
\usepackage{hyperref}
\usepackage{graphicx}
\usepackage{amsmath}
\usepackage{xcolor}
\usepackage{amssymb}
\usepackage{amsfonts}
\usepackage{times}
\usepackage{comment}
\usepackage{marginnote}
\usepackage{bm}
\usepackage{float}

\usepackage[normalem]{ulem}

\newtheorem{theorem}{Theorem}[section]

\newtheorem{remark}{Remark}[section]

\newcommand{\bal}{\begin{aligned}}
\newcommand{\eal}{\end{aligned}}
\def\beq#1{\begin{equation}\label{#1}}
\newcommand{\eeq}{\end{equation}}

\def\mcO{{\mathcal O}}

\def\dexpm1{{2}}

\newcommand{\mcG}{\mathcal{G}}
\newcommand{\mcF}{\mathcal{F}}
\newcommand{\mcI}{\mathcal{I}}

\newcommand{\mcN}{\mathcal{N}}
\newcommand{\mcL}{\mathcal{L}}

\newcommand{\mcX}{\mathcal{X}}

\newcommand{\mbR}{\mathbb{R}}

\newcommand{\mbRn}{{\mathbb{R}^n}}

\def \xb{\mathbf{x}}

\def \zb{\mathbf{z}}

\def \yb{\mathbf{y}}

\renewcommand{\vec}[1]{\mathbf{#1}}
\newcommand{\dd}{\,d}

\allowdisplaybreaks

\newcommand{\abs}[1]            {\left|#1\right|}
\newcommand{\norm}[1]           {\left|\!\left|#1\right|\!\right|}
\newcommand{\opnorm}[1]         {\left|\!\left|\!\left|#1\right|\!\right|\!\right|}
\usepackage{layouts}

\begin{document}
\begin{frontmatter}

\title{An asymptotically compatible coupling formulation for\\ nonlocal interface problems with jumps}

\author{Christian Glusa}
\address{Center for Computing Research,
         Sandia National Laboratories,
         Albuquerque, NM}
         
\author{Marta D'Elia}
\address{Data Science and Computing Group,
         Sandia National Laboratories,
         Livermore, CA}

\author{Giacomo Capodaglio}
\address{Computational Physics and Methods Group,
         Los Alamos National Laboratory,
         Los Alamos, NM}
         
\author{Max Gunzburger}
\address{Department of Scientific Computing,
         Florida State University,
         Tallahassee, FL}
         
\author{Pavel B.~Bochev}
\address{Center for Computing Research,
         Sandia National Laboratories,
         Albuquerque, NM}

\begin{abstract}
We introduce a mathematically rigorous formulation for a nonlocal interface problem with jumps and propose an asymptotically compatible finite element discretization for the weak form of the interface problem. After proving the well-posedness of the weak form, we demonstrate that solutions to the nonlocal interface problem converge to the corresponding local counterpart when the nonlocal data are appropriately prescribed. Several numerical tests in one and two dimensions show the applicability of our technique, its numerical convergence to exact nonlocal solutions, its convergence to the local limit when the horizons vanish, and its robustness with respect to the patch test. 
\end{abstract}

\begin{keyword}
Nonlocal equations, Interface problems, Imperfect interfaces, Finite element discretizations, Coupling, Asymptotic behavior of solutions
\end{keyword}

\end{frontmatter}
\section{Introduction}
Nonlocal models have become popular alternatives to model phenomena where classical partial differential equations fail to be descriptive. These phenomena include the presence of discontinuities in the solution (such as fractures in continuum mechanics \cite{silling2000reformulation}), the appearance of anomalous diffusion effects (such as super- and sub-diffusion in subsurface transport and turbulence \cite{suzuki2021fractional}), and the presence of heterogeneities at the small scales that affect the global behavior of a system at the macro scale (such as the effects of micro-scale heterogeneities in wave propagation \cite{You2021AAAI}). As a result, several scientific and engineering fields have benefited from the use of nonlocal equations. These include, but are not limited to, fracture mechanics \cite{Ha2011,silling2000reformulation}, anomalous subsurface transport \cite{Benson2000,Gulian2021,Schumer2003,Schumer2001,xu2022machine}, image processing \cite{Buades2010,DElia2021Imaging,Gilboa2007,you2021NKN}, stochastic processes \cite{Burch2014,DElia2017,Meerschaert2012,MeKl00}, turbulence \cite{Akhavan2021LES,DiLeoni2021,pang2020npinns,Pang2019fPINNs}, and homogenization \cite{You2020Regression,you2022data}.

The distinguishing feature of nonlocal operators is their integral nature which allows for capturing long-range interactions and reduces the regularity constraints on the solutions. In this paper we limit our attention to spatial nonlocal elliptic operators, specifically to the nonlocal Laplacian operator. In its most general form, it is defined as \cite{DElia2021Unified}
\begin{equation}\label{eq:nl-laplacian}
\mcL u(\xb):= 2 \int_\mbRn (u(\yb)-u(\xb)) \gamma(\xb,\yb) d\yb    
\end{equation}
where the kernel function $\gamma$, usually symmetric and nonnegative, is compactly supported in the {\it nonlocal neighborhood} $B_\delta(\xb)$, i.e. the ball centered at $\xb$ of radius $\delta$. The latter, often referred to as {\it horizon}, determines the extent of the long-range nonlocal interactions and it is allowed to be infinite. Operators such as \eqref{eq:nl-laplacian} have been extensively studied in the last decade and their analysis is supported by a rigorous nonlocal vector calculus theory that has significantly advanced during that time. We refer the reader to \cite{Du2012,Du2013,Lehoucq-Gunzburger} for seminal nonlocal calculus concepts and to \cite{DElia2021Unified,d2021connections} for a more general, unified theory. 

While significant progress has been made in the context of modeling and simulation of nonlocal operators, several mathematical and computational challenges still remain, hindering the use of nonlocal equations in scientific and engineering applications. In this paper we focus on the specific problem of treating physical interfaces that may arise, e.g., in the presence of material discontinuities. While the analysis and discretization of interface problems is well-established for PDEs \cite{cao2010coupled,chen2011parallel,hansbo2002unfitted}, only a few works have addressed the same problem for nonlocal equations and even fewer papers have proposed rigorous interface formulations. In \cite{Alali2015,seleson2013interface} a nonlocal interface formulation was proposed; however, these works do not provide rigorous analysis of either the well-posedness of the formulation or the convergence to the local, PDE limit, as the horizon vanishes. Furthermore, they do not treat the relevant case of solution or flux jumps occurring at the interface. The first rigorous formulation of a nonlocal interface problem can be found in \cite{Capodaglio2020} where the authors propose a formulation based on an energy principle. Here, in the absence of interface jumps in the solution or in the flux, conditions for well-posedness and convergence to the local limit are rigorously derived and illustrated by numerical examples. Related works that treat nonlocal interface problems deal with fractional models \cite{d2021fractional} and domain-decomposition methods for nonlocal equations \cite{aksoylu2011variational,capodaglio2020general,xu2021feti}.

In this paper, building on \cite{Capodaglio2020}, we propose a modified formulation that treats the so-called {\it imperfect interface} problem (see \cite{JAVILI201476} for the local counterpart). Specifically, we introduce the first mathematically rigorous formulation of an interface problem with jumps and an asymptotically compatible discretization. Our major contributions are as follows.
\begin{itemize}
    \item We propose a novel formulation of a nonlocal interface problem with variable horizon and with jumps and we prove that, under certain assumptions on the kernels, it is well-posed. We also show that in the absence of jumps, when the kernel is uniformly defined over the entire domain and across the interface, and when the horizon is constant, our formulation is equivalent to that of a single-domain problem. 
    \item We propose a finite element discretization and show that the $h$-convergence behavior solely depends on the choice of finite element space. In fact, our numerical illustrations in one and two dimensions show optimal quadratic convergence in the $L^2$ norm when using piecewise linear finite element spaces. 
    \item We prove that the solutions of the proposed formulation converge to the solutions of the corresponding local problem when the data of the nonlocal interface problem are properly prescribed. Specifically, given a reference local interface problem we propose nonlocal jump conditions on the solution and the flux that guarantee convergence as the horizon vanishes. We also justify the numerically observed rate of convergence in the $H^1$ seminorm. 
    \item We present several one- and two-dimensional numerical tests that illustrate the applicability of the proposed formulation and its robustness. In particular, we report the results of the so-called patch test that show that the method is consistent in the case of linear solutions. This is a desirable property for every coupling method.
\end{itemize}

\paragraph{Outline of the paper}
In Section \ref{sec:nonlocal-interface} we introduce the proposed strong and weak forms of the nonlocal interface problem and prove the well-posedness of the latter for certain classes of kernels. We provide numerical illustrations of the applicability of this technique in one and two dimensions and show that the interface problem discretized with the finite element method converges to its continuous counterpart as we refine the mesh. We also confirm the robustness of our method by providing results for a linear patch test. In Section \ref{sec:local-limit} we first introduce the weak form of a local interface problem with jumps and show that our formulation converges to its local counterpart in the limit of vanishing nonlocality. We also report several numerical illustrations that confirm the convergence of the nonlocal solutions to the local one as well as a confirmation of the expected convergence rate in the $H^1$ seminorm. In Section \ref{sec:conclusion} we summarize our contributions and provide future research plans.

\section{A nonlocal interface model}\label{sec:nonlocal-interface}
In this section we introduce the strong and weak forms of the nonlocal interface problem and provide a well-posedness result for the latter. 

Although we start our description from the strong form, as it is common in the context of PDE models, we point out that our formulation is derived by mimicking the weak form of a local interface problem with jumps in the solution and in the flux at the interface. Such a local formulation is reported in Section \ref{sec:local-limit}, where we study the convergence of solutions of the nonlocal interface problem to their local counterpart.

\subsection{Strong form}

\begin{figure}[t]
\centering
\includegraphics[width=0.7\textwidth]{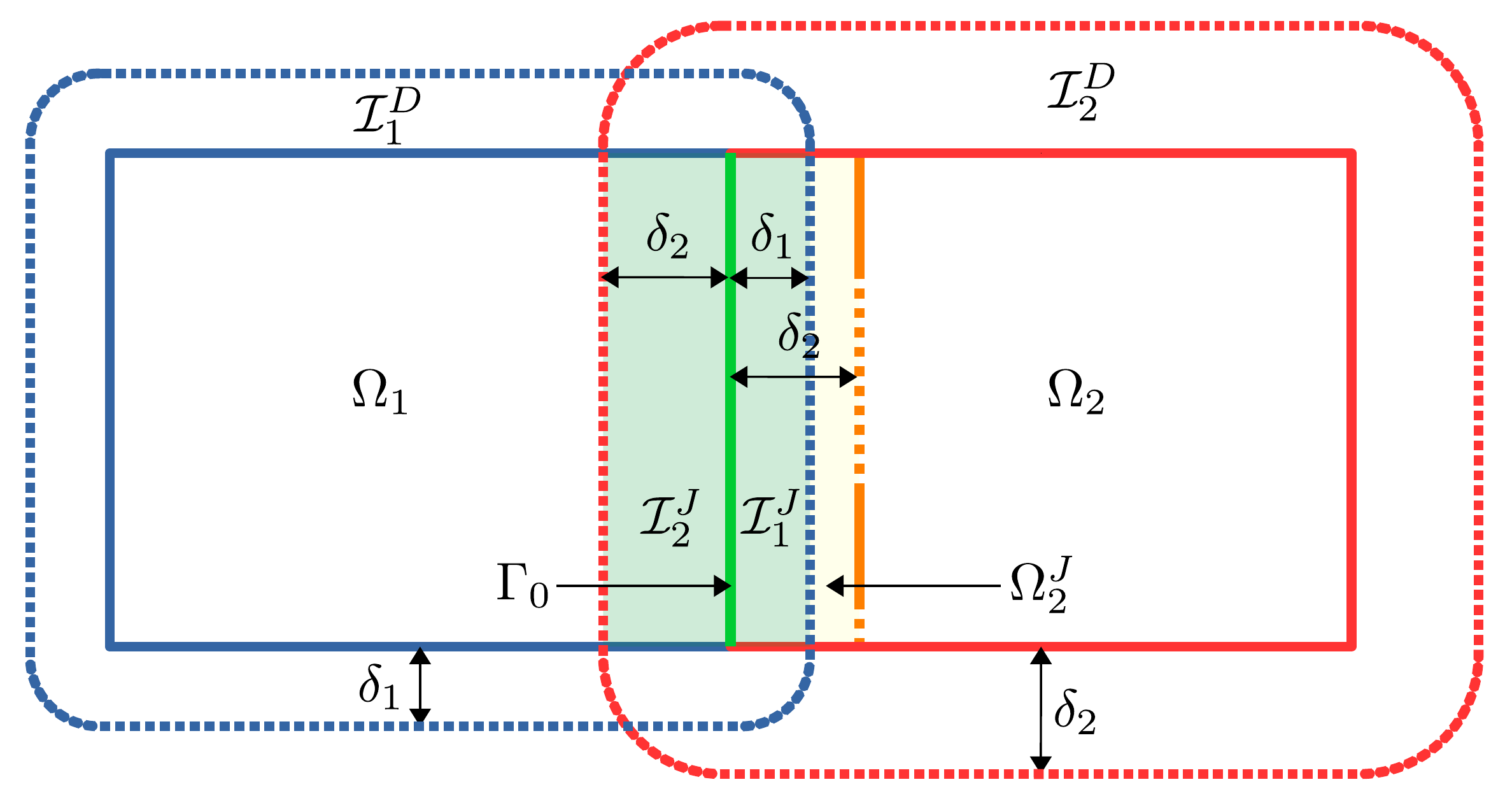}
\\[-0.3em]
\caption{Illustration of a two-dimensional configuration of a nonlocal interface problem with different horizons. The nonlocal interface $\Gamma$ corresponds to $\mathcal I_1^J\cup\mathcal I_2^J$.}
\label{fig:domains}
\end{figure}

Let \(\Omega_{i}\subset\mathbb{R}^{n}\), with \(n=1,2,3\) and $i=1,2$, denote open domains such that \(\Omega_{1}\cap\Omega_{2}=\emptyset\) and \(\partial\Omega_{1}\cap\partial\Omega_{2}=\Gamma_{0}\neq\emptyset\). We refer to $\Gamma_0$ as the {\it local} interface. For a two-dimensional illustration, we refer to the configuration shown in Figure \ref{fig:domains}. 
Let \(\delta_{i}>0\) be two \emph{interaction horizons} and let $|\cdot|$ denote the $\ell^2$ norm; we define the so-called interaction domains as
\begin{align*}
  \mcI_{i}:=\{\yb\in\mathbb{R}^{n}\setminus\Omega_{i} \text{ such that } |\xb-\yb|<\delta_{i} \text{ for some } \xb\in\Omega_{i}\}.
\end{align*}
We assume without loss of generality that \(\delta_{1}\leq \delta_{2}\).
In order to identify the nonlocal interface we partition each \(\mcI_{i}\) into two non-overlapping parts:
\begin{align*}
  \mcI_{i}^{J} &:= \mcI_{i}\cap \Omega_{j}, \\
  \mcI_{i}^{D} &:= \mcI_{i}\setminus\mcI_i^J,
\end{align*}
where \(\mcI_{i}^{D}\) is the portion of interaction domain where nonlocal ``boundary conditions''\footnote{We will refer to the nonlocal counterpart of boundary conditions as volume constraints, introduced later on in the paper.} for the sub-problem on domain \(i\) are available, whereas \(\mcI_{i}^{J}\) is the portion of interaction domain that overlaps with the adjacent domain and where jump conditions (if any) are to be applied. As such, the latter domains form the {\it nonlocal} interface \(\Gamma:=\mcI_{1}^{J}\cup\mcI_{2}^{J}\cup\Gamma_0\).
Here, and in what follows, we will use the index \(j:=3-i\) to denote the subdomain neighboring \(i\).

For scalar functions $u_{i}:\Omega_i\cup\mcI_i\to\mbR$, \(i=1,2\), the nonlocal Laplacian is defined as
\begin{equation}\label{eq:laplacian}
    \mathcal{L}_{i}u_{i}(\xb):= 2\int_{\Omega_{i}\cup\mcI_{i}} (u_{i}(\yb)-u_{i}(\xb)) \gamma_{i}(\xb,\yb;\delta_{i}) d\yb,  \quad \xb\in\Omega_{i},
\end{equation}
where \(\gamma_{i}\) are nonnegative, symmetric kernel functions such that\footnote{Although the Euclidean norm is the most standard choice, we note that the theory developed in this work is not tied to the specific norm used to define the nonlocal neighborhood; for example the $\ell^\infty$ norm is also a choice, see \cite{xu2021machine,xu2021feti} for examples where this neighborhood is used.} \(\gamma_{i}(\xb,\yb;\delta_{i})=0\) for \(|\xb-\yb|\geq\delta_{i}\). 
We also define
\begin{align*}
  \Omega_i^J:=\{\yb\in\Omega_{i}\setminus\mcI_{j}^{J} \text{ such that } |\xb-\yb|<\delta_{i} \text{ for some } \xb \in \mcI_{i}^{J}\}.
\end{align*}
Since we assumed that \(\delta_{1}\leq \delta_{2}\), we have that \(\Omega_1^J\) is always empty. Moreover, in the specific case of $\delta_1=\delta_2$, both $\Omega_i^J$ are empty.
For the same functions, we introduce the associated {\it interface-flux} operator that is piecewise defined as
\begin{equation}\label{eq:flux}
\begin{aligned}
  \mcF (u_{1},u_{2})(\xb) &:=
      \displaystyle 2\int_{\Omega_{i}^{J}} (u_{i}(\xb)-u_{i}(\yb)) \gamma_{i}^{J}(\xb,\yb;\delta_{i}) d\yb & \\[4mm]
    &\quad\displaystyle + 2\int_{\mcI_{j}^{J}} (u_{j}(\xb)-u_{j}(\yb)) \left(\gamma_{j}^{J}(\xb,\yb;\delta_{j}) - \gamma_{j}(\xb,\yb;\delta_{j})\right) d\yb, & \xb\in \mcI_{i}^{J}, \\
  \mcF(u_{1},u_{2})(\xb) &:= \displaystyle 2\int_{\mcI_{2}^{J}} (u_{2}(\xb)-u_{2}(\yb)) \left(\gamma_{2}^{J}(\xb,\yb;\delta_{2}) - \gamma_{2}(\xb,\yb;\delta_{2})\right) d\yb, & \xb\in \Omega_2^J,
\end{aligned}
\end{equation}
with yet to be defined symmetric kernels \(\gamma_{i}^{J}\), \(i=1,2\), such that \(\gamma_{i}^J(\xb,\yb;\delta_{i})=0\) for \(|\xb-\yb|\geq\delta_{i}\). Note that in the original nonlocal vector calculus literature \cite{Du2012}, the flux operator, i.e. the operator associated with Neumann or flux conditions in a single domain problem, is expressed as a composition of the nonlocal interaction operator $\mcN$ and the nonlocal gradient operator $\mcG$; for any point $\xb$ in the interaction domain and for a globally defined kernel function $\gamma$, the resulting operator reads
$$
\mcN(\mcG u)(\xb) = -\displaystyle \int_{\Omega\cup\mcI} (u(\yb)-u(\xb)) \gamma(\xb,\yb)\,d\yb.
$$

The more complex structure of the interface flux in \eqref{eq:flux} stems from arguments related to the energy of the coupled problem. In fact, as we show in the following paragraphs, in the absence of jumps through the interface our formulation recovers the formulation of a single-domain nonlocal problem. More generally, in the presence of jumps, the formulation that we propose corresponds to a single domain problem with a piecewise defined kernel and additional forcing terms appropriately prescribed around the nonlocal interface. This important observation is the key to proving the well-posedness and limiting behavior of the nonlocal interface problem.

As mentioned above, the following strong form of the nonlocal interface problem is posed in such a way that the corresponding weak form mimics its local counterpart. For $i=1,2$, we introduce the forcing terms $\zeta_i:\Omega_i\to\mbR$, the Dirichlet data $\kappa_i:\mcI_i^D\to\mbR$, the solution jump $\mu:\Gamma\to\mbR$, and the flux jump $\nu:\Gamma\cup\Omega_2^J\to\mbR$, and we seek for $u_i:\Omega_i\cup\mcI_i\to\mbR$ such that the following equations hold.
\begin{equation}\label{eq:strong-form}
\begin{aligned}
-\mathcal{L}_{i}u_{i} &= \zeta_{i}(\xb) & \xb&\in\Omega_{i}, \\
  u_{i}(\xb) &= \kappa_{i}(\xb) &  \xb&\in\mcI_{i}^{D},\\
  u_{2}(\xb)-u_{1}(\xb)&=\mu(\xb) &  \xb&\in\Gamma, \\
  \mathcal{F}(u_{1},u_{2})(\xb)&=\nu(\xb) & \;\;\; \xb&\in\Gamma\cup\Omega_2^J.
\end{aligned}
\end{equation}
Here, the second condition is the nonlocal counterpart of a Dirichlet boundary condition; differently from the local case, it is prescribed on a volume and it is referred to as a Dirichlet volume constraint. Neumann constraints have also been considered in the literature \cite{d2020physically,d2021prescription}; however, in this work we only consider the Dirichlet case as the prescription of Neumann conditions is not germane to the paper. We note that the last two conditions in \eqref{eq:strong-form} couple the equations for $u_1$ and $u_2$ and represent jump conditions on the solutions and fluxes; we refer to them as {\it nonlocal interface conditions}. We study the well-posedness of this coupled system in its weak form, introduced in the following section.

\subsection{Weak form}
We derive the weak form of \eqref{eq:strong-form} by multiplying the first row by a test function, $v_i$, suitably chosen to be consistent with the volume constraints on $\mcI_i^D$ and the interface conditions on $\Gamma$. 
For $i=1,2$, let \(v_{i}=0\) on \(\mcI_{i}^{D}\) and \(v_{1}=v_{2}\) in \(\Gamma\). For simplicity, but without loss of generality, in this section we consider the homogeneous case\footnote{The non-homogeneous case can be treated by lifting arguments, as it is often done in the context of non-homogeneous variational problems, see, e.g., \cite{DElia2014DistControl}.} for which $u_i=0$ in $\mcI_i^D$ and $u_1=u_2$ in $\Gamma$, i.e. $\mu=0$.
By adding the first equation in \eqref{eq:strong-form} for $i=1,2$, multiplying by $v_i$, integrating over $\Omega_i$, and using the nonlocal flux condition, we obtain
\begin{equation}\label{eq:weak1}
\begin{aligned}
  &\sum_{i=1}^{2} \left\{\int_{\Omega_{i}}\zeta_{i}v_{i}d\xb + \int_{\mcI_{i}^{J}}\nu v_{i}d\xb\right\} + \int_{\Omega_2^J}\nu v_{2}d\xb \\
  =& \sum_{i=1}^{2} \left\{ 2 \int_{\Omega_{i}} \int_{\Omega_{i}\cup \mcI_{i}} (u_{i}(\xb)-u_{i}(\yb)) v_{i}(\xb)\gamma_{i}(\xb,\yb;\delta_{i}) d\yb d\xb \right. \\
    & \qquad + 2 \int_{\Omega_{i}^{J}} \int_{\mcI_{i}^{J}} (u_{i}(\xb)-u_{i}(\yb)) v_{i}(\xb) \left(\gamma_{i}^{J}(\xb,\yb;\delta_{i})-\gamma_{i}(\xb,\yb;\delta_{i})\right) d\yb d\xb \\
    & \qquad + 2 \int_{\mcI_{i}^{J}} \int_{\Omega_{i}^{J}} (u_{i}(\xb)-u_{i}(\yb)) v_{i}(\xb)\gamma_{i}^{J}(\xb,\yb;\delta_{i}) d\yb d\xb \\
    & \qquad \left.+ 2 \int_{\mcI_{j}^{J}} \int_{\mcI_{i}^{J}} (u_{i}(\xb)-u_{i}(\yb)) v_{i}(\xb) \left(\gamma_{i}^{J}(\xb,\yb;\delta_{i}) - \gamma_{i}(\xb,\yb;\delta_{i})\right) d\yb d\xb \right\}.
  \end{aligned}
\end{equation}
Here, we have used that \(\Omega_{1}^{J}=\emptyset\).
Now, since we have that \(\Omega_i^J\subset\Omega_{i}\setminus\mcI_{j}^{J}\) and that \(\gamma_{i}(\xb,\yb;\delta_{i})=0\) and \(\gamma_{i}^{J}(\xb,\yb;\delta_{i})=0\) for \(\xb\in\Omega_{i}\setminus (\mcI_{j}^{J}\cup\Omega_i^J)\), \(\yb\in\mcI_{i}^{J}\), we can expand the domain of integration of the second term from \(\Omega_{i}^{J}\times\mcI_{i}^{J}\) to \(\left(\Omega_{i}\setminus\mcI_{j}^{J}\right)\times\mcI_{i}^{J}\).
Similarly, we can expand the third term from \(\mcI_{i}^{J}\times\Omega_{i}^{J}\) to \(\mcI_{i}^{J}\times\Omega_{i}\).
With the purpose of simplifying the expression above, we collect the domains of integration after cancellation for the kernels \(\gamma_{i}\) and \(\gamma_{i}^{J}\) separately.
\(\gamma_{i}\) is integrated over
\begin{align*}
  \left[\Omega_{i}\times(\Omega_{i}\cup\mcI_{i})\right] \setminus \left[(\Omega_{i}\setminus\mcI_{j}^{J})\times\mcI_{i}^{J}\right] \setminus \left[\mcI_{j}^{J}\times\mcI_{i}^{J}\right] &= \Omega_{i} \times (\Omega_{i} \cup \mcI_{i}^{D}),
\end{align*}
whereas \(\gamma_{i}^{J}\) is integrated over
\begin{align*}
  \left[(\Omega_{i}\setminus\mathcal{I}_{j}^{J})\times\mcI_{i}^{J}\right] \cup \left[\mcI_{i}^{J}\times\Omega_{i}\right] \cup \left[\mcI_{j}^{J}\times\mcI_{i}^{J}\right] &= \left[\Omega_{i}\times\mcI_{i}^{J}\right] \cup \left[\mcI_{i}^{J}\times\Omega_{i}\right].
\end{align*}
Thus, expression \eqref{eq:weak1} becomes
\begin{equation}\label{eq:weak2}
  \begin{aligned}
    &\sum_{i=1}^{2}\left\{ \int_{\Omega_{i}}\zeta_{i}v_{i}\,d\xb + \int_{\mcI_{i}^{J}}\nu v_{i}\,d\xb\right\} + \int_{\Omega_2^J}\!\!\!\nu v_{2} \,d\xb\\
    =& \sum_{i=1}^{2} \left\{ 2 \iint_{\Omega_{i}\times(\Omega_{i}\cup \mcI_{i}^{D})} (u_{i}(\xb)-u_{i}(\yb)) v_{i}(\xb)\gamma_{i}(\xb,\yb;\delta_{i}) d\yb d\xb \right. \\
    & \qquad \left. + 2 \iint_{\left[\Omega_{i}\times\mcI_{i}^{J}\right] \cup \left[\mcI_{i}^{J}\times\Omega_{i}\right]}  (u_{i}(\xb)-u_{i}(\yb)) v_{i}(\xb) \gamma_{i}^{J}(\xb,\yb;\delta_{i}) d\yb d\xb \right\}.
  \end{aligned}
\end{equation}
To rewrite the above expression in a more compact form, we introduce the functional
\begin{equation}\label{eq:functional}
      r(v_1,v_2)=\sum_{i=1}^{2}\left\{ \int_{\Omega_{i}}\zeta_{i}v_{i}\,d\xb+ \int_{\mcI_{i}^{J}}\nu v_{i}\,d\xb\right\} + \int_{\Omega_2^J}\!\!\!\nu v_{2}\,d\xb,
\end{equation}
the symmetric kernels\footnote{For numerical implementation purposes, we also report the following systematic way of defining $\widetilde\gamma_i$:
\begin{center}
\begin{tabular*}{0.22\linewidth}{r|ccc}
    \(\xb , \yb\) & \(\Omega_{i}\) & \(\mcI_{i}^{J}\) & \(\mcI_{i}^{D}\) \\ \hline
    \(\Omega_{i}\) &   \(\gamma_{i}\) & \(\gamma_{i}^{J}\) & \(\gamma_{i}\)  \\
    \(\mcI_{i}^{J}\) & \(\gamma_{i}^{J}\) & 0 & 0 \\
    \(\mcI_{i}^{D}\) & \(\gamma_{i}\) & 0 & 0
\end{tabular*}
\end{center}}
\begin{align}\label{eq:subdomainKernels}
\tilde{\gamma}_{i}(\xb,\yb;\delta_{i}) &=
\begin{cases}
    0 & \text{on } \left(\mcI_{i}^{J}\cup\mcI_{i}^{D}\right)^{2},\\
    \gamma_{i}^{J}(\xb,\yb) & \text{on } \left[\Omega_{i}\times\mcI_{i}^{J}\right] \cup \left[\mcI_{i}^{J}\times\Omega_{i}\right],\\
    \gamma_{i}(\xb,\yb) & \text{on } (\Omega_{i} \cup \mcI_{i}^{D})^{2}\setminus \left(\mcI_{i}^{D}\right)^{2}
\end{cases}
\end{align}
and the corresponding bilinear form
\begin{equation}\label{eq:bilinear}
     a(u_1,u_2,v_1,v_2;\delta) := \sum\limits_{i=1}^2 \iint_{(\Omega_{i}\cup\mcI_{i})^{2}} (u_{i}(\xb)-u_{i}(\yb)) (v_{i}(\xb)-v_{i}(\yb)) \tilde{\gamma}_{i}(\xb,\yb;\delta_{i})  d\yb d\xb.
\end{equation}
Then, with these definitions, \eqref{eq:weak2} is equivalent to
\begin{equation}\label{eq:weak}
   a(u_1,u_2,v_1,v_2;\delta) = r(v_1,v_2),
\end{equation}
where we used the integration by parts formula \cite{Du2012} and the fact that \(v_{i}=0\) on \(\mcI_{i}^{D}\) and \(v_{1}=v_{2}\) on \(\Gamma\).
We observe that the interactions between the two subdomains \(\Omega_{i}\) are given by \(\gamma_{i}^{J}\) in \eqref{eq:subdomainKernels}.
This explains the particular choice of the flux condition.
Further simplification can be obtained by defining a global solution, $u$, a global test function, $v$, and a global right-hand side $\zeta$. Since $\mu=0$, we set
\begin{align}\label{eq:global-func}
u &=
\begin{cases}
   u_{1} & \text{in } \Omega_{1}\cup\mcI_{1},\\
   u_{2} & \text{in } \Omega_{2}\cup\mcI_{2},
\end{cases}&
v &=
\begin{cases}
   v_{1} & \text{in } \Omega_{1}\cup\mcI_{1},\\
   v_{2} & \text{in } \Omega_{2}\cup\mcI_{2},
\end{cases}&
\zeta =
\left\{
\begin{array}{ll}
\zeta_1 & \;\;\text{in }\Omega_1,\\[1mm]
\zeta_2 & \;\;\text{in }\Omega_{2}.
\end{array}\right.
\end{align}
Then, we can rewrite the weak formulation \eqref{eq:weak} as
\begin{equation}\label{eq:global-weak}
    \widetilde a(u,v;\delta) = \widetilde r(v)
\end{equation}
where
\begin{equation}
\widetilde r(v) = \int_\Omega \zeta \, v\,d\xb  + \int_{\Gamma\cup\Omega_2^J} \!\!\!\nu \, v\,d\xb,
\end{equation}
\begin{equation}\label{eq:global-a}
\begin{aligned}
    \widetilde a(u,v;\delta)  &= \sum_{i=1}^{2} \iint_{(\Omega_{i}\cup\mcI_{i})^{2}} (u(\xb)-u(\yb)) (v(\xb)-v(\yb)) \tilde{\gamma}_{i}(\xb,\yb;\delta_{i})  d\yb d\xb \\
    &= \iint_{\Omega\cup\mcI}(u(\xb)-u(\yb)) (v(\xb)-v(\yb)) \gamma(\xb,\yb;\delta)  d\yb d\xb,
\end{aligned}
\end{equation}
with $\delta=(\delta_1,\delta_2)$ and the symmetric kernel function\footnote{For numerical implementation purposes, we also report the following systematic way of defining $\gamma$
\begin{center}
  \begin{tabular*}{0.8\linewidth}{r|cc|cc|ccc}
    \(\xb , \yb\)                          & \(\Omega_{1}\setminus\mcI_{2}^{J}\) & \(\mcI_{2}^{J}\) & \(\mcI_{1}^{J}\) & \(\Omega_{2}\setminus\mcI_{1}^{J}\) & \(\mcI_{1}^{D}\setminus \mcI_{2}^{D}\) & \(\mcI_{1}^{D} \cap \mcI_{2}^{D}\) & \(\mcI_{2}^{D}\setminus \mcI_{1}^{D}\) \\ \hline
    \(\Omega_{1}\setminus\mcI_{2}^{J}\)    &  \(\gamma_{1}\) & \(\gamma_{1}\) & \(\gamma_{1}^{J}\) & 0 & \(\gamma_{1}\) & \(\gamma_{1}\) & 0 \\
    \(\mcI_{2}^{J}\)                       &  \(\gamma_{1}\)& \(\gamma_{1}+0\) & \(\gamma_{1}^{J}+\gamma_{2}^{J}\) & \(\gamma_{2}^{J}\) & \(\gamma_{1}\) & \(\gamma_{1}+0\) & 0 \\ \hline
    \(\mcI_{1}^{J}\)                       &  \(\gamma_{1}^{J}\)& \(\gamma_{1}^{J}+\gamma_{2}^{J}\)  & \(0+\gamma_{2}\) & \(\gamma_{2}\) & 0 & \(0+\gamma_{2}\) & \(\gamma_{2}\)\\
    \(\Omega_{2}\setminus\mcI_{1}^{J}\)    & 0 & \(\gamma_{2}^{J}\) & \(\gamma_{2}\) & \(\gamma_{2}\) & 0 & \(\gamma_{2}\) & \(\gamma_{2}\) \\ \hline
    \(\mcI_{1}^{D}\setminus \mcI_{2}^{D}\) & \(\gamma_{1}\) & \(\gamma_{1}\) & 0 & 0 & 0 & 0 & 0\\
    \(\mcI_{1}^{D} \cap \mcI_{2}^{D}\)     & \(\gamma_{1}\) & \(\gamma_{1}+0\) & \(0+\gamma_{2}\) & \(\gamma_{2}\)& 0 & 0 & 0\\
    \(\mcI_{2}^{D}\setminus \mcI_{1}^{D}\) & 0 & 0 & \(\gamma_{2}\) & \(\gamma_{2}\)& 0 & 0 & 0
  \end{tabular*}
\end{center}}
\(\gamma:=\tilde{\gamma}_{1}\chi_{(\Omega_{1}\cup\mcI_{1})^{2}}+\tilde{\gamma}_{2}\chi_{(\Omega_{2}\cup\mcI_{2})^{2}}\). 
We stress that $\widetilde a$ is the bilinear form associated with a single-domain problem over $\Omega\cup\mcI$ for the symmetric kernel $\gamma$. This fact facilitates the subsequent well-posedness and local-convergence analysis as it allows us to use well-established nonlocal calculus tools. 

\subsection{Well-posedness of the interface problem}

We analyze the well-posedness of the weak form \eqref{eq:global-weak}. Thanks to the global reformulation as a single-domain problem, we use basic tools of the nonlocal vector calculus theory. We define the {\it energy semi-norm} as 
\begin{equation}\label{eq:opnorm}
\opnorm{v}^2_\delta :=
\iint_{(\Omega\cup\mcI)^{2}}  (v(\xb)-v(\yb))^2 \gamma(\xb,\yb;\delta)  d\yb d\xb = \widetilde a(v,v;\delta)
\end{equation}
and the corresponding {\it energy space} as
\begin{equation}\label{eq:Espace}
    V_\delta:=\{ v\in L^2(\Omega\cup\mcI): \opnorm{v}_\delta<\infty, \; v_i=0\; {\rm in} \; \mcI_i^D \;
    {\rm and} \; v_1=v_2 \; {\rm in} \; \Gamma\}.
\end{equation}

Some considerations are in order. Based on the nonlocal calculus theory \cite{Du2012}, for certain choices of kernel functions it is possible to show that the space defined in \eqref{eq:Espace} corresponds to well-known Sobolev spaces and that the energy norms are equivalent to the corresponding Sobolev norms. In fact, when the kernel is of fractional type, i.e.
\begin{equation*}
\gamma(\xb,\yb;\delta)\propto |\xb-\yb|^{-(n+2s)}\mcX(|\xb-\yb|<\bar\delta)
\qquad \text{(type 1)}
\end{equation*}
with $\bar\delta=\max_i\delta_i$ and $s\in(0,1)$, then $\opnorm{\cdot}_{\delta}$ is equivalent to $\|\cdot\|_{H^s(\Omega\cup\mcI)}$, whereas when the kernel function is square integrable, i.e.,
\begin{equation*}
    \int_{\Omega\cup\mcI}\gamma(\xb,\yb;\delta)^2 d\yb < \infty, 
    \qquad \text{(type 2)} 
\end{equation*}
the norm $\opnorm{\cdot}_{\delta}$ is equivalent to $\|\cdot\|_{L^2(\Omega\cup\mcI)}$. In this work we consider both fractional (type 1) and square-integrable (type 2) kernels.

Furthermore, an immediate consequence of these equivalence properties is the fact that the semi-norm \eqref{eq:opnorm} is indeed a norm in $V_\delta$ and the space is complete in that metric \cite{Du2012}. 

Based on the definitions above, we let $V'_\delta$ be the dual space of $V_\delta$ and state the weak form of the interface problem as follows: 

{\it let $\zeta\in V_\delta'$ and $\nu\in V_\delta'$, we seek for $u\in V_\delta$ such that, for all $v\in V_\delta$, \eqref{eq:global-weak} holds.}

Note that, while the definition of the flux $\mcF$ in \eqref{eq:flux} is such that $\nu$ is not defined over $\Omega$, an extension to zero over the entire domain allows us to write $\nu\in V'_\delta$. 
We can finally prove the following well-posedness result for the kernel classes introduced above.

\begin{theorem}
For kernels of type 1 and 2 and for $(\zeta_1,\zeta_2)\in V_\delta'$ and $\nu\in V_\delta'$, there exists a unique $u\in V_\delta$ such that, for all $v\in V_\delta$, \eqref{eq:global-weak} holds. 
\end{theorem}

\noindent{\it Proof.} Since $\widetilde a(u,u;\delta)=\opnorm{u}_\delta^2$, the bilinear form is coercive in $V_\delta$.
For the same reason, the bilinear form is also continuous in $V_\delta$. We next find a bound for the right-hand side of \eqref{eq:global-weak}. 
We have
\begin{equation*}
\begin{aligned}
|\widetilde r(v)| = \left|\int_{\Omega} \zeta\, v d\xb + \int_{\Gamma\cup\Omega_2^J} \!\!\!\nu\, vd\xb\right| \leq
C \left(\opnorm{\zeta}_{V'_\delta}+ \opnorm{\nu}_{V'_\delta}\right)\,\opnorm{v}_\delta
\end{aligned}
\end{equation*}
where $\opnorm{\cdot}_{V'_\delta}$ represents the norm in the dual space $V'_\delta$ defined in the standard manner via duality pairing and $C$ is a positive constant independent on the function $v$. The inequality above implies that $\widetilde r$ is a continuous functional in $V_\delta$. This fact and the properties of $\widetilde a(\cdot,\cdot;\delta)$ are sufficient conditions for the well-posedness of the weak form \eqref{eq:global-weak}. $\square$

\subsection{A specific choice of the jump kernel function}\label{sec:gamma}
We discuss a possible choice of the kernel function that allows us to simplify the structure of the flux term and, consequently, of the interface problem. We also introduce further assumptions that guarantee the convergence to a local interface problem as the nonlocality vanishes. 

First, we discuss how to define the ``jump'' kernel \(\gamma_{i}^{J}\). Let \(\gamma_{i}^{J}=\frac{1}{2}\gamma_{i}\) on \(\mcI_{1}^{J}\times\mcI_{2}^{J} \cup \mcI_{2}^{J}\times\mcI_{1}^{J}\) and \(\gamma_{i}^{J}=\gamma_{i}\) on \((\Omega_{i}\setminus\mcI_{j}^{J})\times\mcI_{i}^{J} \cup \mcI_{i}^{J}\times(\Omega_{i}\setminus\mcI_{j}^{J})\). Then, the interface-flux operator reduces to
\begin{equation}\label{eq:flux-special-case}
\begin{aligned}
  \mcF(u_{1},u_{2})(\xb)
  &= \displaystyle 2\int_{\Omega_i^J} (u_{i}(\xb)-u_{i}(\yb)) \gamma_{i}(\xb,\yb;\delta_{i}) d\yb & \\[4mm]
  &+\displaystyle \int_{\mcI_{j}^{J}} (u_{i}(\xb)-u_{i}(\yb)) \gamma_{i}(\xb,\yb;\delta_{i}) d\yb \\[4mm]
  &-\displaystyle \int_{\mcI_{j}^{J}} (u_{j}(\xb)-u_{j}(\yb)) \gamma_{j}(\xb,\yb;\delta_{j}) d\yb, & \xb\in \mcI_{i}^{J} \\
  \mcF(u_{1},u_{2})(\xb) &= \displaystyle 0, & \xb\in\Omega_2^J.
\end{aligned}
\end{equation}

With this choice and for \(i=1\) the first integral vanishes since \(\Omega^J_1\) is empty; also, the last two integrals can be interpreted as the difference of two nonlocal fluxes from \(\mcI_{j}^{J}\) to \(\mcI_{i}^{J}\). In particular, this implies that the solution jump \(\mu\) and the flux jump \(\nu\) both will be supported on \(\Gamma\). 

For the remainder of the paper, we will adopt the definition \eqref{eq:flux-special-case} and we will further assume that, for $i=1,2$, the kernel functions $\gamma_i$ are such that for \(\xb\in\Omega_{i}\)
  \begin{equation}\label{eq:gamma-scaling}
    \int_{(\xb+\zb)\in(\Omega_i\cup\mcI_i)} \!\!z_{k} z_{l} \gamma_i(\xb,\xb+\zb) \,d\zb
    = \begin{cases}
      1 & \text{ if } k=l, \\
      0 & \text{ otherwise.}
    \end{cases}
  \end{equation}
This normalization condition will play an important role in the convergence behavior of the nonlocal solution as $\delta_i\to 0$, for $i=1,2$. We also point out that all kernels used in our numerical tests satisfy \eqref{eq:gamma-scaling}.
Furthermore, we note that, with the specific choice of \(\gamma_{i}^{J}\), the above condition also holds for \(\tilde{\gamma}_{i}\), but only for \(\xb\in\Omega_{i}\setminus\mcI_{j}^{J}\).

\subsection{Numerical illustrations and $h$-convergence}\label{sec:numerics-illustrations}
In this section we report one- and two-dimensional illustrations of numerical solutions to the interface problem \eqref{eq:strong-form} in the presence of jumps in the solutions and fluxes on the nonlocal interface. To test the robustness of our implementation and to show that the interface formulation does not alter the convergence behavior of numerical solutions of nonlocal problems, we perform convergence tests with respect to the discretization parameter that show that numerical solutions of \eqref{eq:strong-form} feature an optimal convergence rate. We also perform the so-called patch test that further confirms the robustness of the proposed formulation by showing that linear solutions are reproduced at machine-precision accuracy when a piecewise linear finite element space is utilized.

\paragraph{Discretization} 
In all our tests we utilize a finite element (FE) discretization with piecewise linear basis functions.
For more details on the discrete variational formulation and its convergence properties, we refer the reader to \cite{DElia-ACTA-2020}.
The well-posedness of the FE formulation is inherited from the well-posedness of the continuous problem \eqref{eq:global-weak}.
Note that for standard (i.e. in the absence of interfaces) nonlocal problems, our choice of FE basis yields a quadratic convergence, in the $L^2$ norm, of the discretization error with respect to the FE grid size, which we refer to as $h$.
From now on, we refer to the discretized nonlocal solution as $u_h$ and to the convergence with respect to the mesh size as $h$-convergence. 

\begin{figure}
  \centering
  \includegraphics[width=0.4\textwidth]{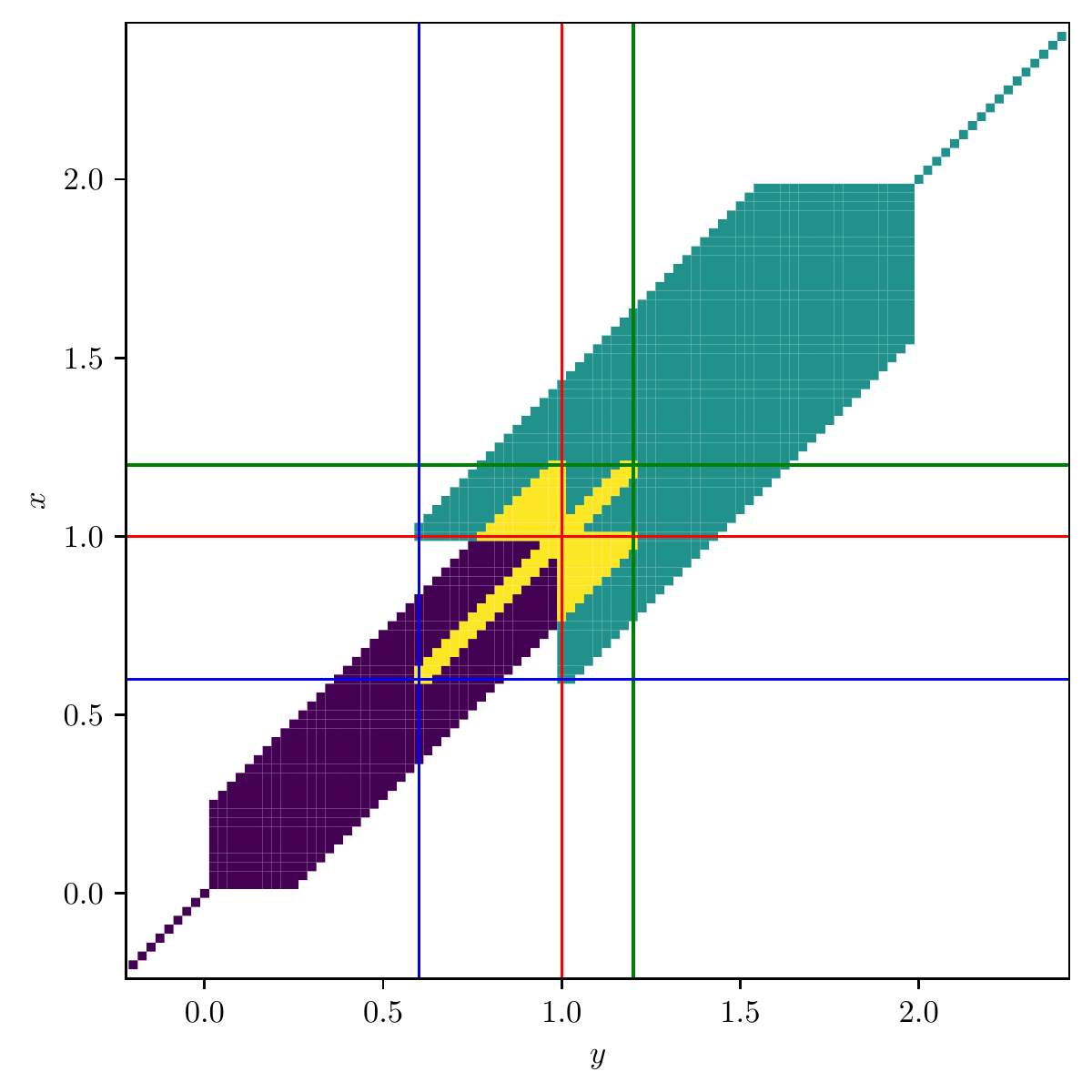}
  \caption{
    Sparsity pattern of \(a(\cdot,\cdot;\delta)\) for the 1D problem over a finite element space of mesh size \(h=2.5\times10^{-2}\) with \(\delta_{1}=0.2\) and \(\delta_{2}=0.4\).
    The local interface \(\Gamma_{0}\) is given by horizontal and vertical red lines, the boundaries of the nonlocal interface \(\Gamma\) by blue and green lines.
    Matrix entries depending only on \(\gamma_{1}\) are shown in dark violet, only on \(\gamma_{2}\) in turquoise, and matrix entries that depend on both in yellow.}
  \label{fig:sparsity-NL-NL}
\end{figure}

To visualize the structure of the discretization matrix associated with the interface problem, in Figure~\ref{fig:sparsity-NL-NL}, we display the sparsity pattern associated with \(\widetilde a(\cdot,\cdot;\delta)\) in a one-dimensional setting. Here, we consider $\Omega_{1}= (0,1)$, $\Omega_{2}=(1,2)$, $\delta_{1}= 0.2$ and $\delta_{2}=0.4$ and a mesh of size \(h=2.5\times10^{-2}\). The local interface \(\Gamma_{0}\) is highlighted by the horizontal and vertical red lines, and the boundaries of the nonlocal interface \(\Gamma\) highlighted by the blue and green lines for $\Omega_1$ and $\Omega_2$ respectively. We report the entries depending only on \(\gamma_{1}\) in dark violet, the ones depending only on \(\gamma_{2}\) in turquoise, and the ones that depend on both in yellow. This figure clearly highlights what parts of the domain interact with each other both within the same domain and across the interface.

In what follows, we use a direct solver to compute numerical solutions.
It should be pointed out though that the structure of the interface problem naturally lends itself to the development of domain decomposition methods.

\paragraph{$h$-Convergence with type 1 kernels in 1D}
We first illustrate the proposed formulation in a one-dimensional setting. Let $\Omega_{1}= (0,1)$, $\Omega_{2}=(1,2)$, $\delta_{1}= 0.2$ and $\delta_{2}=0.4$. We consider a global kernel $\gamma$, defined as in Section \ref{sec:gamma}, of fractional type (i.e. type 1) with the following choice of $\gamma_i$:
\begin{equation}\label{eq:fractional-kernel}
  \gamma_{i}^{\text{1D},F}(\xb,\yb)=C_{\text{1D},s_i,\delta_i}|\xb-\yb|^{-1-2s_i}\,\mcX(|\xb-\yb|<\delta_i),
\end{equation}
with $C_{\text{1D},s_i,\delta_i}= \frac{2-2s_i}{2}\delta_i^{2s_i-2}$, $s_{1}= 0.2$, and $s_{2}=0.4$.
We select the data of the interface problem by using the following manufactured solutions:
\begin{align}
  u_{1}= \sin(\pi\xb), \quad u_{2}= 1-\sin(\pi\xb), \label{eq:exactSolution}
\end{align}
for which the flux jump is computed from \eqref{eq:flux-special-case}, and the source terms, the solution jump, and the Dirichlet volume constraints are obtained by substitution.

\begin{figure}[t]
  \centering
  \includegraphics[width=0.32\textwidth]{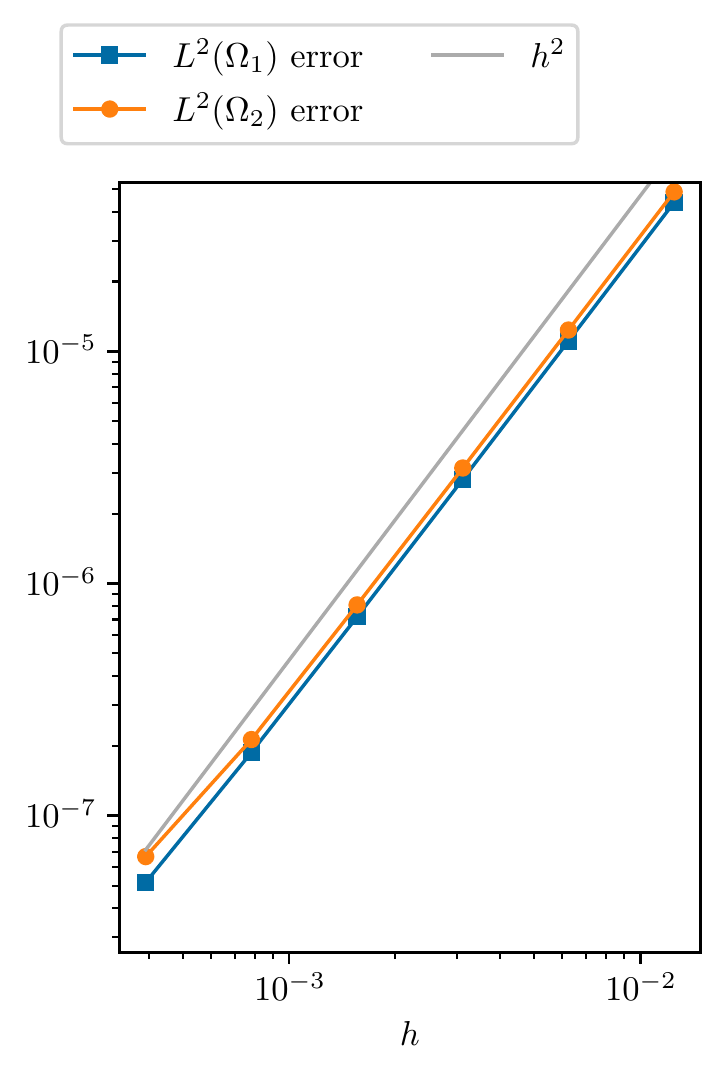}
  \includegraphics[width=0.32\textwidth]{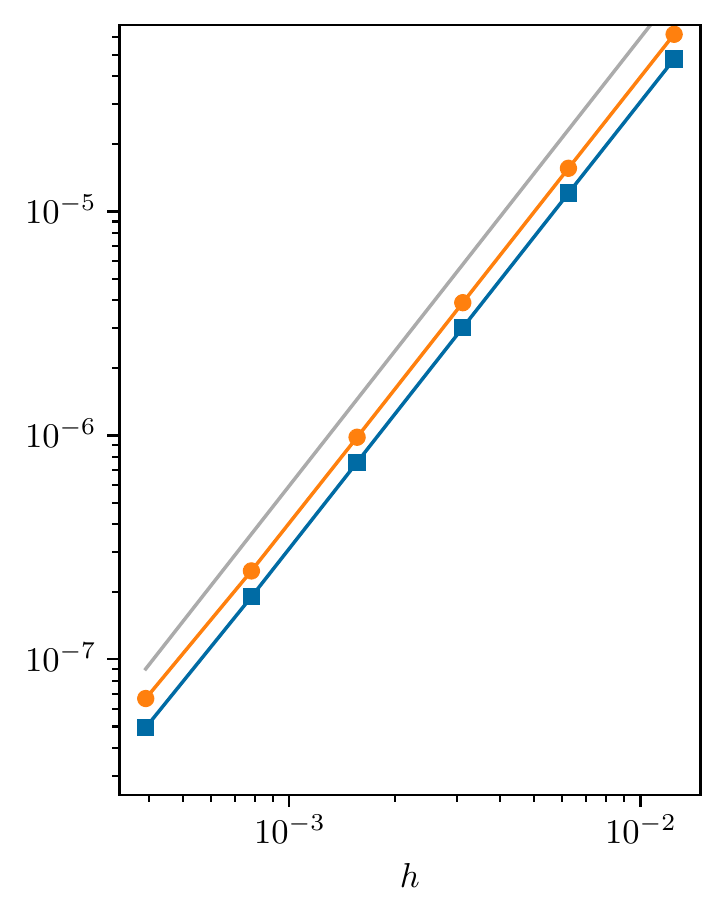}
  \includegraphics[width=0.32\textwidth]{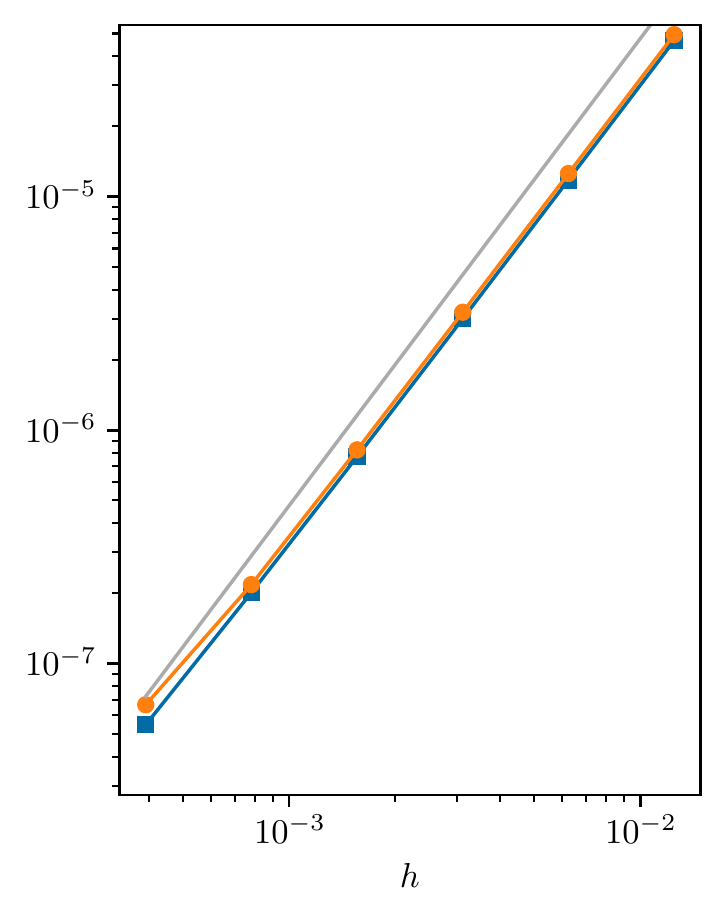}
  \caption{Convergence with respect to the mesh size \(h\).
    \(\delta_{1}=0.2\), \(\delta_{2}=0.4\).
    \emph{Left:} fractional kernels on both subdomains of orders \(s_{1}=0.2\) and \(s_{2}=0.4\).
    \emph{Center:} constant kernels on both subdomains.
    \emph{Right:} constant kernel on the left subdomain, fractional kernel of order \(s_{2}=0.4\) on the right subdomain.}
  \label{fig:hConvergence}
\end{figure}

In Figure~\ref{fig:hConvergence} (left) we show the results of numerical tests for decreasing values of the mesh size $h$. Our $h$-convergence results indicate that the convergence is optimal, i.e. the discretization error converges quadratically in the $L^2$ norm for a smooth exact solution.

\paragraph{$h$-Convergence with type 2 kernels in 1D}
With the same settings and manufactured solution of the previous example, we consider a global integrable kernel $\gamma$, defined as in Section \ref{sec:gamma}, of type 2 with the following choice of $\gamma_i$:
\begin{equation}\label{eq:constant-kernel}
  \gamma_{i}^{\text{1D},C}(\xb,\yb)=C_{\text{1D},\delta_{i}}\mcX(|\xb-\yb|<\delta_{i}),
\end{equation}
with $C_{\text{1D},\delta_{i}}= \frac32\delta_{i}^{-3}$. The forcing and jump terms and the volume constraints are computed, again, by direct substitution of $u_1$ and $u_2$.

In Figure~\ref{fig:hConvergence} (center) we show the results of numerical tests for decreasing values of the mesh size $h$. As in the previous example, our results indicate that the convergence is optimal, i.e. the discretization error converges quadratically in the $L^2$ norm.

\paragraph{$h$-Convergence with kernels of different types in 1D}

We consider different types of kernels on opposite sides of the interface, i.e. we use the constant kernel \(\gamma_{1}^{\text{1D},C}\) in $\Omega_1$, and the fractional kernel \(\gamma_{2}^{\text{1D},F}\) in $\Omega_2$.
Similarly to the previous examples, the source terms and the flux jumps are prescribed so that the analytic solution is again given by \eqref{eq:exactSolution}.

In Figure~\ref{fig:hConvergence} (right), our convergence results confirm that the rate is not affected by the choice of kernel function, i.e. the converge is still optimal.

\paragraph{A linear patch test in 1D} 
We test the consistency of the proposed formulation by showing that the interface formulation passes the so-called patch test. The goal of this test is to show that in the absence of forcing terms and jumps at the interface and when appropriate volume constraints are prescribed, when all the kernel functions are identical, the nonlocal solution is a linear function (as it would happen in a single-domain problem). Furthermore, when the linear solution belongs to the discretization space, as is the case for the piecewise linear discretization utilized in this experiment, the nonlocal interface solution is expected to be accurate up to machine precision. 

\begin{figure}[t]
  \centering
  \includegraphics[width=0.4\textwidth]{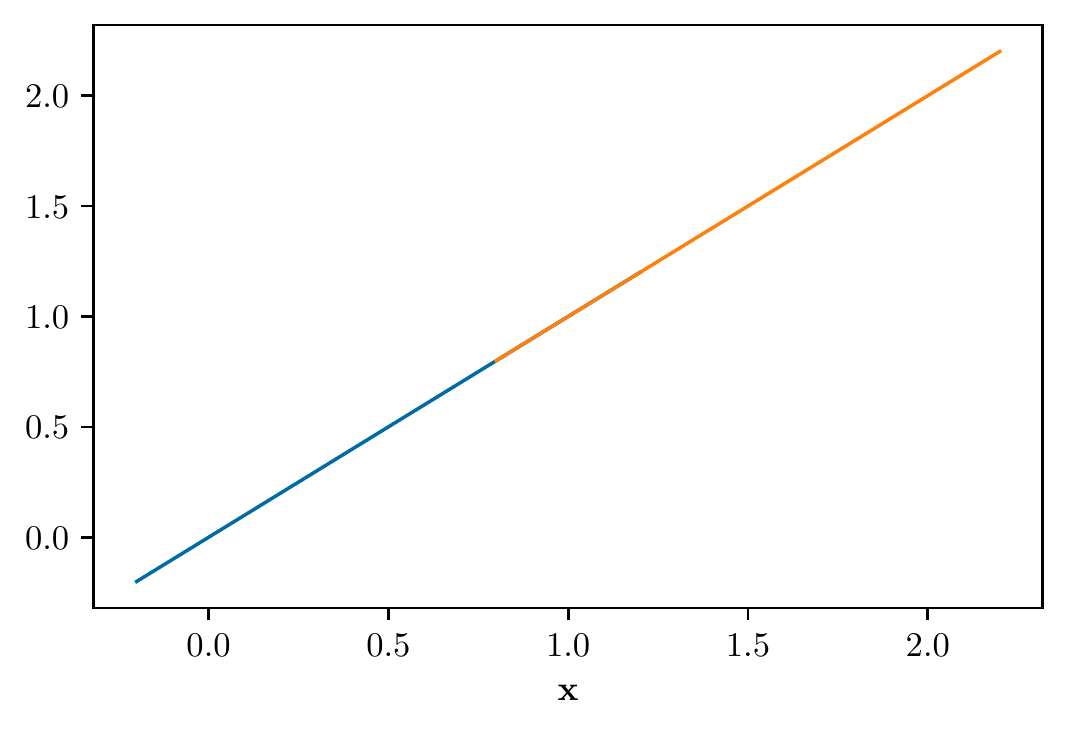}
  \includegraphics[width=0.4\textwidth]{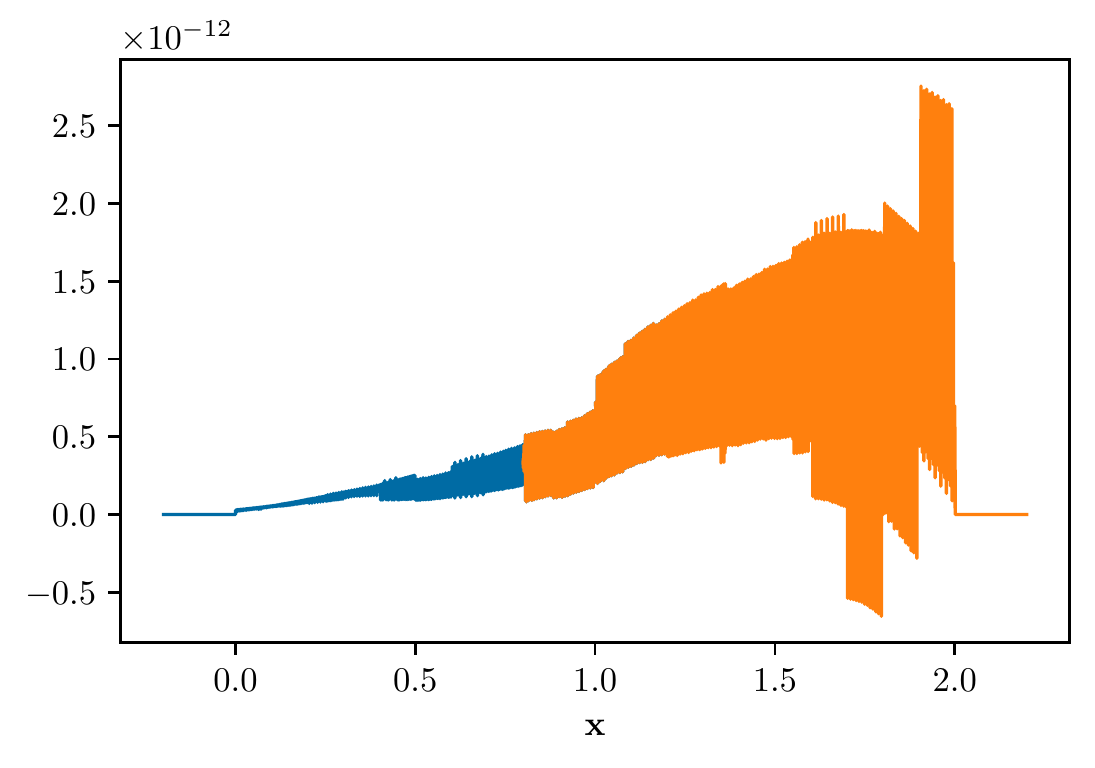}
  \caption{
    Left: Solutions of the patch test for the 1D problem over a finite element space of mesh size \(h \sim 2\times 10^{-4}\) with \(\delta_{1}=\delta_{2}=0.2\) and constant kernels.
    Right: Difference between solution and linear of the patch test.
    }
  \label{fig:patch-constant}
\end{figure}

In Figure~\ref{fig:patch-constant} we report the (linear) nonlocal interface solution obtained for the kernel functions $\gamma_i^{\text{1D},C}$, with $\delta_i=0.2$ for $i=1,2$ and in the one-dimensional setting utilized in the previous experiments. Here, $\mu=0$, $\nu=0$, and $\kappa_i=x$, for $i=1,2$. On the left, we observe the expected linear behavior; on the right, the pointwise error with respect to the exact solution $u=x$ confirms that the nonlocal solution is machine-precision accurate.
A similar behavior, not shown for the sake of brevity, is observed when using fractional kernels with identical orders \(s_{i}\) and horizons \(\delta_{i}\), \(i=1,2\).

\begin{figure}[t!]
  \centering
  \includegraphics[width=0.4\textwidth]{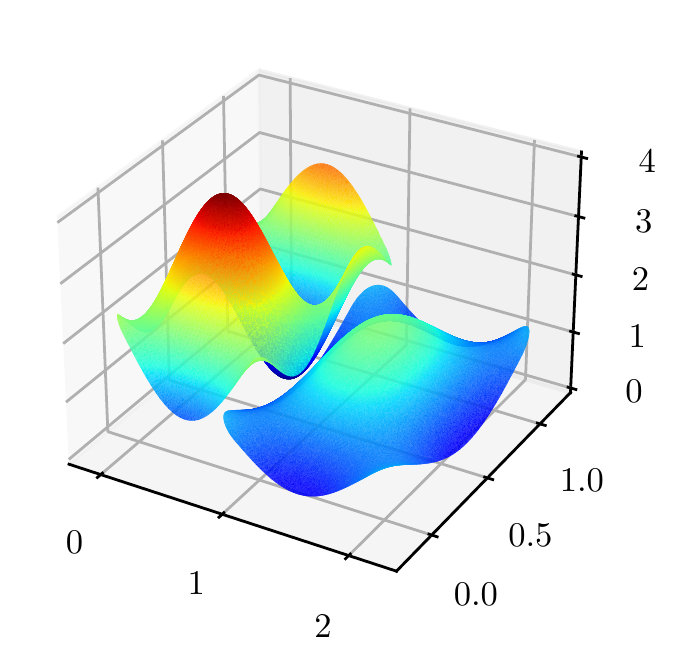}
  \caption{Numerical solution of the nonlocal two-dimensional problem with constant kernels and \(\delta_{1}=0.05\), \(\delta_{2}=0.1\).}
  \label{fig:solPlot2d}
\end{figure}

\paragraph{$h$-Convergence illustrations in 2D}
In this two-dimensional experiment, we let $\Omega_{1}= (0,1)^{2}$ and $\Omega_{2}=(1,2)\times(0,1)$. We set the horizons to $\delta_1=0.1$ and $\delta_2=0.2$. We consider kernels of type 1 and 2 according to the definition of Section \ref{sec:gamma}. For type 1, we choose \(\gamma_{i}\) as
\begin{equation}\label{eq:fractional-kernel2d}
  \gamma_{i}^{\text{2D},F}(\xb,\yb)=C_{\text{2D},s_i,\delta_i}|\xb-\yb|^{-2-2s_i}\,\mcX(|\xb-\yb|<\delta_i),
\end{equation}
with $C_{\text{2D},s_i,\delta_i}= \frac{2-2s_i}{\pi}\delta_i^{2s_i-2}$ for \(s_{1}=0.2\) and \(s_{2}=0.4\).
For type 2, we choose \(\gamma_{i}\) as
\begin{equation}\label{eq:constant-kernel2d}
  \gamma_{i}^{\text{2D},C}(\xb,\yb)=C_{\text{2D},\delta_{i}}\mcX(|\xb-\yb|<\delta_{i}),
\end{equation}
with $C_{\text{2D},\delta_{i}}= \frac{4}{\pi}\delta_{i}^{-4}$.
We select the data of the interface problem by using the following manufactured solutions:
\begin{align}
  u_{1}= 2+2\sin(\pi x_{1})\sin(2\pi x_{2}), \quad u_{2}= 1-\sin(\pi x_{1})\sin(\pi x_{2}), \label{eq:exactSolution2d}
\end{align}
for which, again, the jump and the source terms and the volume constraints are computed using the solutions \(u_1\) and $u_2$.
We report the corresponding solution in Figure \ref{fig:solPlot2d}.

\begin{figure}[t!]
  \centering
  \includegraphics[width=0.32\textwidth]{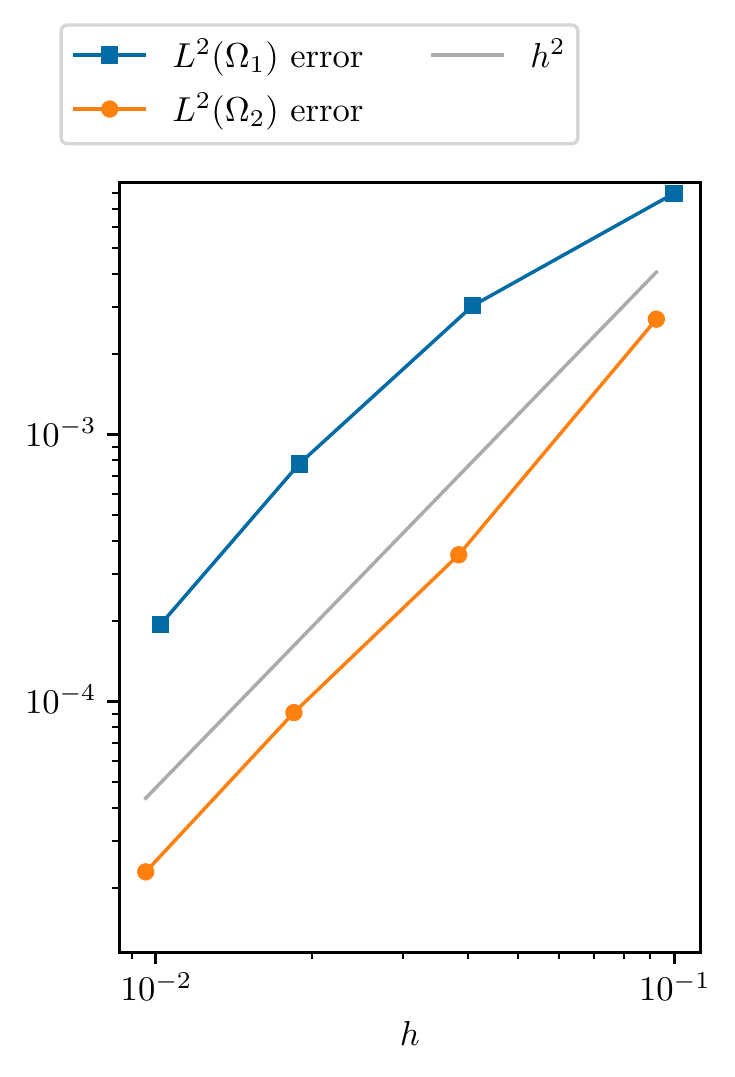}
  \includegraphics[width=0.32\textwidth]{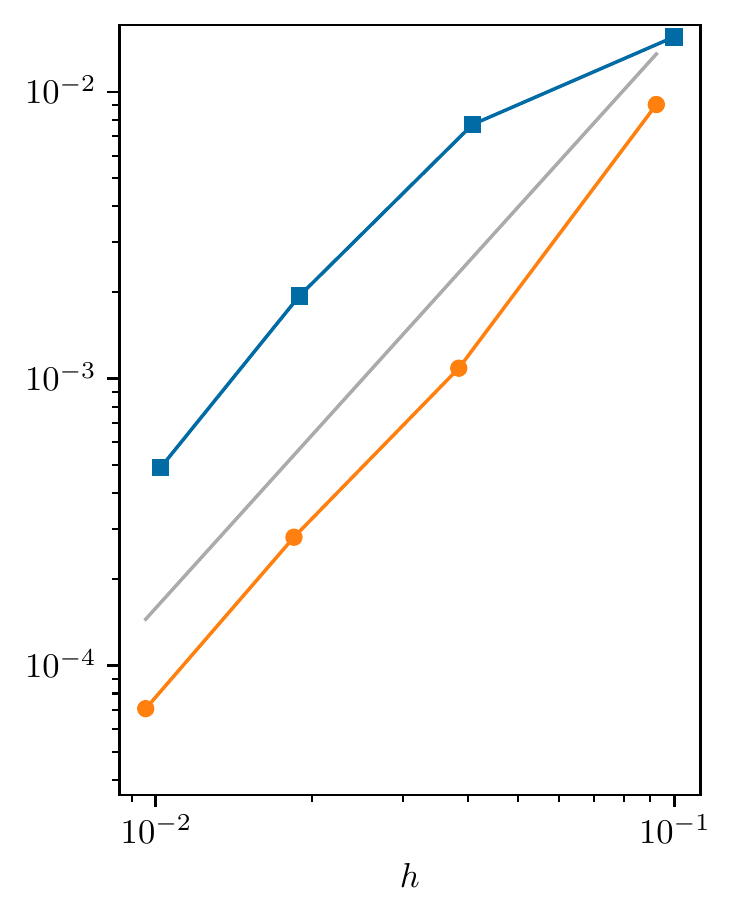}
  \includegraphics[width=0.32\textwidth]{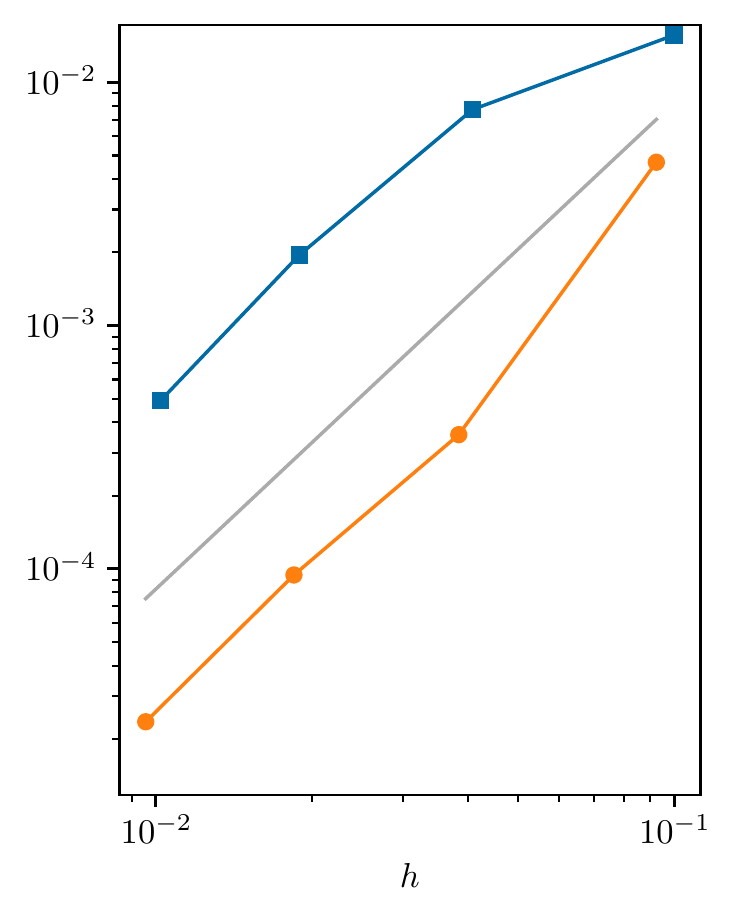}
  \caption{Convergence with respect to the mesh size \(h\).
    \(\delta_{1}=0.1\), \(\delta_{2}=0.2\).
    \emph{Left:} fractional kernels on both subdomains of orders \(s_{1}=0.2\) and \(s_{2}=0.4\).
    \emph{Center:} constant kernels on both subdomains.
    \emph{Right:} constant kernel on the left subdomain, fractional kernel of order \(s_{2}=0.4\) on the right subdomain.
  }
  \label{fig:hConvergence2d}
\end{figure}

In Figure~\ref{fig:hConvergence2d} we display the $h$-convergence to the analytic solution.
It can be observed that as in the previous one-dimensional test cases, the optimal, quadratic $h$-convergence rate is achieved.

\section{Convergence to the local limit}\label{sec:local-limit}

We start this section by recalling the formulation of a local interface problem with jumps in the solutions and in the fluxes at the interface. Next, we present the main result of this section where we describe how to select the nonlocal jump conditions that guarantee the convergence to the local interface problem with prescribed jumps. In the same result we also motivate that solutions to the nonlocal interface problem converge to solutions to the local interface problem in the local energy norm with rate $\mcO(\sqrt{\delta_i})$ with respect to the horizons. Finally, we report several numerical tests that illustrate our theoretical findings. 

\begin{figure}[t]
\centering
\includegraphics[width=0.5\textwidth]{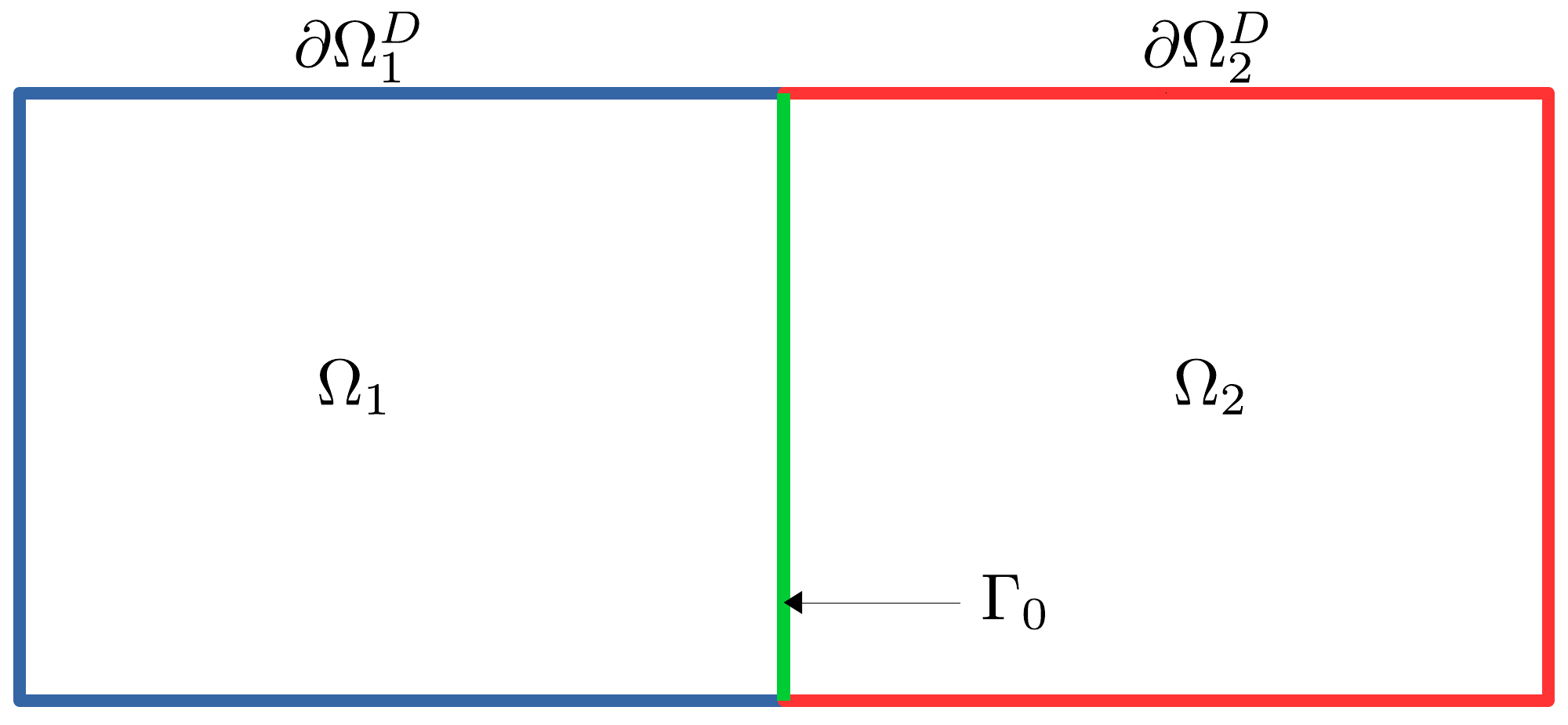}
\\[-0.3em]
\caption{Illustration of a two-dimensional configuration of a local interface.}
\label{fig:local-domains}
\end{figure}

\subsection{A local interface problem with jumps}

We refer to the configuration in Figure \ref{fig:local-domains}, where a two-dimensional illustration is shown. Here, $\Omega_i$, $i=1,2$ are defined as in the previous section; we refer to their ``physical'' boundaries, where Dirichlet boundary conditions are available, as $\partial\Omega_i^D=\partial\Omega_i\setminus\Gamma_0$, with $\Gamma_0=\partial\Omega_1\cap\partial\Omega_2$.

As the local interface problem with jumps is well-known in the literature \cite{JAVILI201476}, we only report its weak form. Let the local bilinear form $\widetilde a(\cdot,\cdot;0):(H^1(\Omega_1)\times H^1(\Omega_2))\times(H^1_0(\Omega_1)\times H^1_0(\Omega_2))\to\mbR$ be defined as 
\begin{equation}\label{eq:local-bilinear}
  a(u_1^0,u_2^0,v_1,v_2;0)=\sum_{i=1}^{2}\int_{\Omega_{i}}\nabla u^0_{i} \cdot \nabla v_{i}\, d\xb,
\end{equation}
where the superscript 0 denotes the local solution. Note that, although we use the same notation used for the nonlocal bilinear form, \eqref{eq:local-bilinear} does not correspond to \eqref{eq:bilinear} evaluated at $\delta=0$, as it is clear from \eqref{eq:local-bilinear}. Then, for $i=1,2$, given $g_i\in H^{1/2}(\Gamma_0)$, $f\in H^{-1}(\Omega)$, and $s,\,h\in H^{1/2}(\Gamma_0)$, we seek $(u_1^0,u_2^0)\in H^1(\Omega_1)\times H^1(\Omega_2)$ such that, for all $(v_1,v_2)\in H^1_0(\Omega_1)\times H^1_0(\Omega_2)$,
\begin{equation}\label{eq:weak-form-local}
a(u_1^0,u_2^0,v_1,v_2;0) = 
\int_\Omega f\, v \,d\xb+ 
\int_{\Gamma_{0}} s \,v \,d\xb,
\end{equation}
subject to
\begin{equation}
\begin{aligned}
u^0_{2}-u^0_{1} & =m &  \;\;\xb&\in\Gamma_{0}, &\\
v_{2}-v_{1} & =0  &\;\;\xb&\in\Gamma_{0}, &\\
u^0_{i} & =g_{i} &  \;\;\xb&\in\partial\Omega_{i}^{D} & \;\;i=1,2, \\
v_{i} & =0 & \;\;\xb&\in\partial\Omega_{i}^{D} & \;\;i=1,2.
\end{aligned}
\end{equation}
In the absence of solution jumps at the interface, i.e. $m=0$, we can introduce the same notation used in the previous section and rewrite the local solutions and tests functions in a more compact way as follows:
\begin{align}\label{eq:global-local-func}
u^0 &=
\begin{cases}
   u^0_{1} & \text{in } \Omega_{1}\cup\Gamma_0,\\
   u^0_{2} & \text{in } \Omega_{2},
\end{cases}&
v &=
\begin{cases}
   v_{1} & \text{in } \Omega_{1}\cup\Gamma_0,\\
   v_{2} & \text{in } \Omega_{2}.
\end{cases}
\end{align}
With this notation, to be consistent with the previous section, we rewrite the bilinear form as 
\begin{equation*}
\widetilde a(u^0,v;0) = a(u_1^0,u_2^0,v_1,v_2;0),
\end{equation*}
so that expression \eqref{eq:weak-form-local} becomes 
\begin{equation}\label{eq:local-globa-weak}
\widetilde a(u^0,v;0) = \int_\Omega f\, v \,d\xb+ 
\int_{\Gamma_{0}} s \,v \,d\xb.
\end{equation}
As mentioned in the previous section, this problem has the same structure as \eqref{eq:global-weak}. The results in the next section confirm that our choice is consistent in the limit of vanishing nonlocal interactions, i.e. solutions to \eqref{eq:global-weak} converge to solutions to \eqref{eq:weak-form-local} as $\delta\to (0,0)$ in the metric induced by the local energy norm. For all functions $v\in H^1(\Omega_1)\times H^1(\Omega_2)$, we define the local energy norm as follows
\begin{equation}\label{eq:local-energy-norm}
\opnorm{v}^2_0 =  \sum_{i=1}^{2}\int_{\Omega_{i}}\left|\nabla v_{i}\right|^2 d\xb.
\end{equation}
Note that $\opnorm{\cdot}_0$ defines a norm equivalent to the norm induced by the local bilinear form.

\subsection{A choice of jump conditions that guarantees convergence as the horizons vanish}

In this section we show how to choose the jump conditions for the nonlocal solutions and fluxes so that the convergence to the local solutions is guaranteed. Also in this case, to avoid technicalities in our analysis, we consider the case of $m=0$ and $\mu=0$ as well as \(g_{i}=0\) and \(\kappa_{i}=0\) and provide conditions for such case in Remark \ref{rem:mu} at the end of this section. Specifically, for $f_i = f|_{\Omega_i}$, given a local problem with data $f_i$ and $s$, we propose choices of $\zeta_i$ and $\nu$, so that solutions to \eqref{eq:global-weak} converge to solutions to \eqref{eq:weak-form-local}. We also provide a result that supports the observed 0.5 convergence rate with respect to $\delta_i$ in the $\opnorm{\cdot}_0$ metric. These statements are summarized in the following theorem.

We define the extension of the local solution \(u^{0}\) to \(\Omega\cup\mcI\) as follows:
\begin{align*}
  Eu^{0}(\xb) :=
  \begin{cases}
    u^{0}(\xb) & \text{if } \xb \in \Omega_{1}\cup\Omega_{2}\cup\Gamma, \\
    0 & \text{else}.
  \end{cases}
\end{align*}

\begin{theorem}\label{thm:delta-convergence}
  Let $\gamma$ be a kernel of type 1 or 2. For $i=1,2$, let $f_i\in H^{-1}(\Omega_i)$ and $s\in H^{1/2}(\Gamma_0)$, and denote by $u^0\in H^1(\Omega_1)\times H^1(\Omega_2)$ the solution to \eqref{eq:local-globa-weak}.
  Choose $\zeta_i$ and $\nu$ such that
\begin{equation}\label{eq:compatible-choice}
\begin{aligned}
\zeta_i & = f_i & \qquad\text{on }\Omega_i, \\[1mm]
\nu(\xb) & = \frac{s(\Pi_{\Gamma_0} \xb)}{\delta_{1}+\delta_{2}} & \qquad\xb\in\Gamma,
\end{aligned}    
\end{equation}
where $\Pi_{\Gamma_0}$ is a projection from $\Gamma$ onto $\Gamma_0$.
If the solution of the nonlocal interface problem \eqref{eq:global-weak} with data as in \eqref{eq:compatible-choice} is such that
\begin{itemize}
\item \(u, Eu^{0} \in C^{\alpha}(\Omega_{i}\cup\mcI_{i})\) where \(\alpha>\max\{1/2,s\}\) for type 1 kernels and \(\alpha>1/2\) for type 2 kernels,
\item \(u, u^{0} \in W^{1,\infty}(\Omega_{i})\), and
\item \(u\) and \(u^{0}\) are \(C^{2}(\Omega_{i}\setminus\mcI_{j}^{J})\),
\end{itemize}
then
\begin{equation}\label{eq:loc-limit}
  \lim\limits_{\delta\to(0,0)}\opnorm{u-u^0}_0 = 0.
\end{equation}
\end{theorem}

\noindent{\it Proof.}
We introduce the ``error function'' \(e:=u-u^0\), defined on $\Omega_1\cup\Omega_2\cup\Gamma_{0}$.
Then, for $i=1,2$ and $\xb\in\partial\Omega_{i}^{D}$, we have $e =0$.
Since \(u\) has by assumption enough regularity, the error function \(e\) can be used as a test function for the local problem.

We also define \(Ee:=u-Eu^{0}\).
We have that $Ee=0$ in $\mcI^D_1\cup\mcI^D_2$.
Therefore, \(Ee\) can be used as a test function for the nonlocal problem.

\smallskip
Our goal is to find a bound for the norm $\opnorm{e}_0$.
We have that
\begin{align*}
  \opnorm{e}_0^{2}
  &\leq \widetilde a(e,e;0) \\
  &= \widetilde a(u,e;0) - \widetilde a(u^0,e;0) \\
  &= \widetilde a(u,Ee;\delta) - \widetilde a(u^0,e;0) + \left(\widetilde a(u,e;0)-\widetilde a(u,Ee;\delta)\right) \\
  &= \int_{\Omega_{1}} (\zeta_{1}-f_{1})e\,d\xb + \int_{\Omega_{2}} (\zeta_{2}-f_{2})e \,d\xb \\
  &+ \int_{\Gamma} \nu Ee\,d\xb
  - \int_{\Gamma_{0}} s e \,d\xb+ \left(\widetilde a(u,e;0)-\widetilde a(u,Ee;\delta)\right).
\end{align*}
By choosing \(\zeta_{i}=f_{i}\), according to the compatibility conditions \eqref{eq:compatible-choice}, we can eliminate the first two terms, and hence obtain
\begin{align}
  \opnorm{e}_0^{2}
  &\leq \abs{\int_{\Gamma} \nu Ee\,d\xb
  - \int_{\Gamma_{0}} s e \,d\xb\;} + \abs{\widetilde a(u,e;0)-\widetilde a(u,Ee;\delta)}. \label{eq:errorBoundDelta}
\end{align}
In what follows, we estimate both terms in \eqref{eq:errorBoundDelta} separately.
In order to avoid technicalities, we assume that \(\Gamma_{0}\) is straight.
Write $\xb\in\Gamma$ as \(\xb=\vec{w}+x\vec{n}\) where \(\vec{n}\) is the normal of \(\Gamma_{0}\) and \(\vec{w}\in\Gamma_{0}\).
We choose \(\nu(\xb):=\frac{1}{\delta_{1}+\delta_{2}}s(\Pi_{\Gamma_{0}}\xb)\) and split the integration domain as \(\Gamma=\Gamma_{0}\times((-\delta_{2},0) \cup [0,\delta_{1}))\) and using the definition of the extension operator \(E\),
\begin{align*}
  \int_{\Gamma} \nu(\xb) Ee(\xb) \dd\xb
  &= \frac{1}{\delta_{1}+\delta_{2}}\int_{\Gamma_{0}} s(\vec{w}) \int_{-\delta_{2}}^{\delta_{1}}  Ee(\vec{w}+x\vec{n}) \dd x \dd\vec{w}\\
  &= \frac{1}{\delta_{1}+\delta_{2}}\int_{\Gamma_{0}} s(\vec{w}) \left[\int_{0}^{\delta_{1}} e(\vec{w}+x\vec{n}) \dd x + \int_{-\delta_{2}}^{0} e(\vec{w}+x\vec{n}) \dd x \right] \dd\vec{w}.
\end{align*}
Via Taylor expansion around \(\vec{w}\) up to the first order derivatives, we have that
\begin{align*}
  e(\vec{w}+x\vec{n}) &= e(\vec{w}) + \int_{0}^{x} (\partial_{\vec{n}}e)(\vec{w}+z\vec{n}) \dd z, & \text{for } x>0,\\
  e(\vec{w}+x\vec{n}) &= e(\vec{w}) + \int_{x}^{0} (\partial_{\vec{n}}e)(\vec{w}+z\vec{n}) \dd z, & \text{for } x<0.
\end{align*}
This holds by density of $C^\infty(\Omega_i)$ in $H^1(\Omega_i)$.
Here, we treated separately the cases \(x>0\) and \(x<0\), since \(Ee\) is generally only continuous across \(\Gamma_{0}\). Integration along \(x\) leads to
\begin{align*}
   & \int_{0}^{\delta_{1}} e(\vec{w}+x\vec{n}) \dd x + \int_{-\delta_{2}}^{0} e(\vec{w}+x\vec{n}) \dd x\\
  =& \delta_{1}e(\vec{w}) + \delta_{2}e(\vec{w}) + \int_{0}^{\delta_{1}} \int_{0}^{x} (\partial_{\vec{n}} e)(\vec{w}+z\vec{n}) \dd z \dd x + \int_{-\delta_{2}}^{0} \int_{x}^{0} (\partial_{\vec{n}} e)(\vec{w}+z\vec{n}) \dd z \dd x.
\end{align*}
Consequently,
\begin{align*}
  &\abs{\int_{\Gamma} \nu(\xb) Ee(\xb) \dd\vec{x} - \int_{\Gamma_{0}} s(\vec{w}) e(\vec{w}) \dd\vec{w}} \\
  =& \frac{1}{\delta_{1}+\delta_{2}}\abs{\int_{\Gamma_{0}} s(\vec{w}) \left[\int_{0}^{\delta_{1}} \int_{0}^{x} (\partial_{\vec{n}} e)(\vec{w}+z\vec{n}) \dd z \dd x + \int_{-\delta_{2}}^{0} \int_{x}^{0} (\partial_{\vec{n}}e)(\vec{w}+z\vec{n}) \dd z \dd x \right] \dd\vec{w}}.
\end{align*}
Using the Cauchy-Schwarz inequality, first with respect to \(\vec{w}\) and then with respect to \(x\), we obtain
\begin{align*}
  &(\delta_{1}+\delta_{2})\abs{\int_{\Gamma} \nu(\xb) (Ee)(\xb) \dd\vec{x} - \int_{\Gamma_{0}} s(\vec{w}) e(\vec{w}) \dd\vec{w}} \\
  \leq &\norm{s}_{L^{2}(\Gamma_{0})} \left[\sqrt{\int_{\Gamma_{0}} \left(\int_{0}^{\delta_{1}} \int_{0}^{x} (\partial_{\vec{n}} e)(\vec{w}+z\vec{n}) \dd z \dd x \right)^{2} \dd\vec{w}}\right. \\ &\hspace{1.5cm}\left.+\sqrt{\int_{\Gamma_{0}} \left(\int_{-\delta_{2}}^{0} \int_{x}^{0} (\partial_{\vec{n}}e)(\vec{w}+z\vec{n}) \dd z \dd x\right)^{2} \dd\vec{w}}\right] \\
  \leq & \norm{s}_{L^{2}(\Gamma_{0})} \left[\sqrt{\delta_{1} \int_{\Gamma_{0}} \int_{0}^{\delta_{1}} \left(\int_{0}^{x} (\partial_{\vec{n}} e)(\vec{w}+z\vec{n}) \dd z\right)^{2} \dd x \dd\vec{w}} \right.\\
  &\hspace{1.5cm}\left.+ \sqrt{\delta_{2} \int_{\Gamma_{0}} \int_{-\delta_{2}}^{0} \left(\int_{x}^{0} (\partial_{\vec{n}}e)(\vec{w}+z\vec{n}) \dd z\right)^{2} \dd x \dd\vec{w}}\right],
\end{align*}
where, in the last step, we have used that \(\int_{0}^{\delta_{1}}\dd x=\delta_{1}\) and that \(\int_{-\delta_{2}}^{0}\dd x=\delta_{2}\). Again by the Cauchy-Schwarz inequality, we have that
\begin{align*}
  &\left(\int_{0}^{x} (\partial_{\vec{n}} e)(\vec{w}+z\vec{n}) \dd z\right)^{2}\\
  \leq &\left(\int_{0}^{x} \dd z\right)  \left(\int_{0}^{x} \abs{(\partial_{\vec{n}} e)(\vec{w}+z\vec{n})}^{2} \dd z\right)\\
  \leq& \delta_{1}\int_{0}^{\delta_{1}} \abs{(\partial_{\vec{n}} e)(\vec{w}+z\vec{n})}^{2} \dd z,
\end{align*}
where we have used that \(0\leq x\leq \delta_{1}\) in the first integral.
Using the same argument for the second term, we obtain
\begin{align*}
  &\abs{\int_{\Gamma} \nu Ee- \int_{\Gamma_{0}} s e} \\
  & \leq \frac{1}{\delta_{1}+\delta_{2}} \norm{s}_{L^{2}(\Gamma_{0})} \left[\sqrt{\delta_{1}^{2} \int_{\Gamma_{0}} \int_{0}^{\delta_{1}} \int_{0}^{\delta_{1}} \abs{(\partial_{\vec{n}} e)(\vec{w}+z\vec{n})}^{2} \dd z \dd x\dd\vec{w}} \right.\\
  &\hspace{3cm}\left.+ \sqrt{\delta_{2}^{2} \int_{\Gamma_{0}} \int_{-\delta_{2}}^{0} \int_{-\delta_{2}}^{0} \abs{(\partial_{\vec{n}}e)(\vec{w}+z\vec{n})}^{2} \dd z \dd x \dd\vec{w}}\right] \\
  & = \frac{1}{\delta_{1}+\delta_{2}} \norm{s}_{L^{2}(\Gamma_{0})} \left[\sqrt{\delta_{1}^{3} \int_{\Gamma_{0}} \int_{0}^{\delta_{1}} \abs{(\partial_{\vec{n}} e)(\vec{w}+z\vec{n})}^{2} \dd z \dd\vec{w}} \right.\\
  &\hspace{3cm}\left.+\sqrt{\delta_{2}^{3} \int_{\Gamma_{0}} \int_{-\delta_{2}}^{0} \abs{(\partial_{\vec{n}}e)(\vec{w}+z\vec{n})}^{2} \dd z \dd\vec{w}}\right] \\
  & = \frac{1}{\delta_{1}+\delta_{2}} \norm{s}_{L^{2}(\Gamma_{0})} \left[\delta_{1}^{3/2}\norm{\partial_{\vec{n}}e}_{L^{2}(\mcI^J_2)} + \delta_{2}^{3/2}\norm{\partial_{\vec{n}}e}_{L^{2}(\mcI^J_1)}\right].
\end{align*}
Since \(\norm{\partial_{\vec{n}}e}_{L^{2}(\mcI^J_2)},\norm{\partial_{\vec{n}}e}_{L^{2}(\mcI^J_1)}\leq \opnorm{e}_{0}\), we finally obtain
\begin{align}\label{eq:jump-limit}
  \abs{\int_{\Gamma} \nu(\xb) Ee(\xb) \dd\xb - \int_{\Gamma_{0}} s(\vec{w}) e(\vec{w}) \dd\vec{w}} &\leq (\sqrt{\delta_{1}}+\sqrt{\delta_{2}}) \norm{s}_{L^{2}(\Gamma_{0})}\opnorm{e}_{0}
\end{align}
for the first term of \eqref{eq:errorBoundDelta}.
For the second term of \eqref{eq:errorBoundDelta}, we use standard nonlocal vector calculus arguments. We have
\begin{equation}\label{eq:second-term}
\begin{aligned}
& \;\abs{\widetilde a(u,e;0)-\widetilde a(u,Ee;\delta)} \\
= & \;\abs{\sum\limits_{i=1}^2\left\{\int_{\Omega_i} \nabla u \nabla Ee\, d\xb\right\} - \iint_{(\Omega\cup\mcI)^2} (u(\yb)-u(\xb))(Ee(\yb)-Ee(\xb))\gamma(\xb,\yb)\,d\xb\,d\yb}
\\[2mm]
= &  \;\abs{\sum\limits_{i=1}^2\left\{\int_{\Omega_i} \nabla u \nabla Ee\, d\xb -
\iint_{(\Omega_{i}\cup\mcI_{i})^{2}} (u(\xb)-u(\yb)) (Ee(\xb)-Ee(\yb)) \tilde{\gamma}_{i}(\xb,\yb;\delta_{i})  d\yb d\xb\right\}}.
\end{aligned}
\end{equation}
We partition
\begin{align*}
  \Omega_{i}\cup\mcI_{i} = \Lambda_{i}\cup \left[\left(\Omega_{i}\cup\mcI_{i}\right)\setminus\Lambda_{i}\right],
\end{align*}
where
\begin{align*}
  \Lambda_{i}:=\{\xb \in \Omega_{i} \;|\; \forall \yb\in \mcI_{i}\cup\Gamma: |\xb-\yb|>\delta_{i}\} \subset \Omega_{i}\setminus\mcI_{j}^{J}.
\end{align*}
Consequently, we can write
\begin{align*}
  &\iint_{(\Omega_{i}\cup\mcI_{i})^{2}} (u(\xb)-u(\yb)) (Ee(\xb)-Ee(\yb)) \tilde{\gamma}_{i}(\xb,\yb;\delta_{i})  d\yb d\xb \\
  =&\int_{\Lambda_{i}}\int_{\Omega_{i}\cup\mcI_{i}} (u(\xb)-u(\yb)) (Ee(\xb)-Ee(\yb)) \tilde{\gamma}_{i}(\xb,\yb;\delta_{i})  d\yb d\xb \\
  &+\int_{(\Omega_{i}\cup\mcI_{i})\setminus\Lambda_{i}}\int_{\Omega_{i}\cup\mcI_{i}} (u(\xb)-u(\yb)) (Ee(\xb)-Ee(\yb)) \tilde{\gamma}_{i}(\xb,\yb;\delta_{i})  d\yb d\xb.
\end{align*}
In the first term, we can Taylor expand since \(B_{\delta}(\xb)\subset\Omega_{i}\setminus\mcI_{j}^{J}\) and we have by assumption \(C^{2}\) regularity of \(u\) and \(Ee|_{\Omega_{i}\setminus\mcI_{j}^{J}}=u-u^{0}\):
\begin{align*}
  u(\yb) &= u(\xb) + \nabla u(\xb)\cdot (\yb-\xb) + \mathcal{O}(\delta^{2}) \\
  Ee(\yb) &= Ee(\xb) + \nabla Ee(\xb)\cdot (\yb-\xb) + \mathcal{O}(\delta^{2}),
\end{align*}
and obtain
\begin{align*}
  &\int_{\Lambda_{i}}\int_{\Omega_{i}\cup\mcI_{i}} (u(\xb)-u(\yb)) (Ee(\xb)-Ee(\yb)) \tilde{\gamma}_{i}(\xb,\yb;\delta_{i})  d\yb d\xb \\
  =& \int_{\Lambda_{i}} \nabla u(\xb) \cdot \left[\int_{\Omega_{i}\cup\mcI_{i}} (\xb-\yb)(\xb-\yb)\tilde{\gamma}_{i}(\xb,\yb) d\yb\right] \nabla Ee(\xb) d\xb  + \mathcal{O}(\delta_{i}).
\end{align*}
Moreover, the $(k,l)$ component of the tensor inner integral over $\Omega_i\cup\mcI_i$, is given by
\begin{align*}
  [{\rm int}]_{kl}=\int_{\Omega_i\cup\mcI_i} (\yb-\xb)(\yb-\xb)\tilde\gamma_i(\xb,\yb) \,d\yb
  &= \left[\int_{(\xb+\zb)\in(\Omega\cup\mcI)} \!\!\! z_{k} z_{l} \tilde\gamma_i(\xb,\xb+\zb) \,d\zb\right]_{kl}\!.
\end{align*}
For \(\xb\in\Lambda_{i}\) we have that \(\tilde\gamma_i(\xb,\xb+\cdot)=\gamma_{i}(\xb,\xb+\cdot)\).
Since \(\gamma_{i}\) is symmetric and scaled according to \eqref{eq:gamma-scaling}, we have that
\begin{align*}
  [{\rm int}]_{kl}= \begin{cases}
    1 & \text{ if } k=l, \\
    0 & \text{ otherwise.}
  \end{cases}
\end{align*}
Thus, we have
\begin{equation}\label{eq:part1}
\begin{aligned}
  &\int_{\Lambda_{i}}\int_{\Omega_{i}\cup\mcI_{i}} (u(\xb)-u(\yb)) (Ee(\xb)-Ee(\yb)) \tilde{\gamma}_{i}(\xb,\yb;\delta_{i})  d\yb d\xb \\
  =& \int_{\Lambda_{i}} \nabla u(\xb) \cdot \nabla Ee(\xb) d\xb  + \mathcal{O}(\delta_{i}).
\end{aligned}
\end{equation}
For the second term, we use that by assumption \(u_{i}\in C^{\alpha}(\Omega_{i}\cup\mcI_{i})\) and write
\begin{equation}\label{eq:part2}
\begin{aligned}
  &\left|\int_{(\Omega_{i}\cup\mcI_{i})\setminus\Lambda_{i}}\int_{\Omega_{i}\cup\mcI_{i}} (u(\xb)-u(\yb)) (Ee(\xb)-Ee(\yb)) \tilde{\gamma}_{i}(\xb,\yb;\delta_{i})  d\yb d\xb \right| \\
  \leq& \norm{u}_{C^{\alpha}}\norm{Ee}_{C^{\alpha}} \int_{(\Omega_{i}\cup\mcI_{i})\setminus\Lambda_{i}}\int_{\Omega_{i}\cup\mcI_{i}} \abs{\xb-\yb}^{2\alpha} \tilde{\gamma}_{i}(\xb,\yb;\delta_{i})  d\yb d\xb \\
  \leq& C_i\norm{u}_{C^{\alpha}}\norm{Ee}_{C^{\alpha}} \delta_{i}^{2\alpha-1},
\end{aligned}
\end{equation}
where the constants $C_i$ are independent of \(\delta_{i}\). Thus, by combining \eqref{eq:part1} and \eqref{eq:part2} we have that \eqref{eq:second-term} is equivalent to
\begin{align*}
&\abs{\widetilde a(u,e;0)-\widetilde a(u,Ee;\delta)}\\
=&  \abs{\sum_{i=1}^2\left\{\int_{\Omega_{i}\setminus\Lambda_{i}}\nabla u\cdot\nabla Ee + \mcO(\delta_i) \right\}}\\[2mm]
\leq & \sum_{i=1}^2\left\{C_i \delta_{i} \norm{\nabla u}_{L^{\infty}}\norm{\nabla Ee}_{L^{\infty}} + \mcO(\delta_i)\right\},
\end{align*}
where the last inequality follows from the fact that \(\abs{\Omega_{i}\setminus\Lambda_{i}}=\mathcal{O}(\delta_{i})\) and \(W^{1,\infty}\) boundedness of \(u\) and \(Ee=u-Eu^{0}\). Thus, since we assumed that \(\alpha>1/2\), the following equality implies that the error converges to zero as $\delta\to(0,0)$:
\begin{align*}
  \abs{\widetilde a(u,e;0)-\widetilde a(u,Ee;\delta)} = \mathcal{O}(\delta_{1}^{2\alpha-1}+\delta_{2}^{2\alpha-1}).
\end{align*}
$\square$

\begin{remark}\label{rem:conv-rate}
The estimate \eqref{eq:jump-limit} indicates that the convergence rate in the local energy metric of the error cannot be faster than $\mathcal O(\sqrt{\delta_1}+\sqrt{\delta_2})$. 
For sufficient solution regularity, we expect the convergence to be dominated by the first term in \eqref{eq:errorBoundDelta}, i.e.
\begin{equation}\label{eq:loc-limit-rate}
  \opnorm{u-u^0}_0 \leq c (\sqrt{\delta_1}+\sqrt{\delta_2}),
\end{equation}
where the constant \(c\) depends on the data, but not on \(\delta_{i}\). Indeed, for smooth local solutions,
this estimate seems sharp, as confirmed by the numerical tests reported in the following section.
\end{remark}

\begin{remark}\label{rem:mu}
Note that when the horizons $\delta_i$ vanish, the nonlocal interface $\Gamma$ collapses into $\Gamma_0$. In the presence of solution jumps and non-homogeneous volume conditions, the following condition on $\mu$ and \(\kappa_{i}\) guarantees the convergence of $u$ to $u^0$ for $\delta_i\to 0$, as illustrated in the next section. For $m\in H^{\frac12}(\Gamma_0)$, we set
$$\mu|_{\Gamma_0}  = m,$$
and for \(g_{i}\in H^{1/2}(\partial\Omega^{D}_{i})\) we require
\begin{align*}
  \kappa_{i}|_{\partial\Omega^{D}_{i}} = g_{i}.
\end{align*}
\end{remark}

\subsection{Numerical $\delta$-convergence}

Using the same discretization introduced in the previous section, we test the convergence of the nonlocal solution to the nonlocal interface problem to its local counterpart as the horizons $\delta=(\delta_1,\delta_2)$ go to zero. We refer to this type of convergence as $\delta$-convergence. 

\paragraph{Problem setting in 1D} 
We consider the same one-dimensional setting utilized in the previous section and let $\Omega_{1}= (0,1)$ and $\Omega_{2}=(1,2)$. We first define the local interface problem and then choose the data of the nonlocal interface problem according to Theorem \ref{thm:delta-convergence}. Let the local source terms, the Dirichlet boundary conditions, and the jump terms be defined as
\begin{align*}
  f_{1}&= \pi^{2}\sin(\pi \xb), & f_{2}&= -2\pi^{2}\sin(\pi \xb)\\
  g_{1}&= \sin(\pi\xb), & g_{2}&= 1-2\sin(\pi\xb)\\
  s&= -3\pi, &
  m&=1.
\end{align*}
The exact solution to the local problem is then given by $u_{1}= \sin(\pi\xb)$ and $u_{2}= 1-2\sin(\pi\xb)$. The corresponding nonlocal source terms, Dirichlet volume constraints, and jump terms are then given by
\begin{align}
\zeta_{1}&= \pi^{2}\sin(\pi \xb), & \zeta_{2}&= -2\pi^{2}\sin(\pi \xb), \nonumber\\
  \kappa_{1}&= \sin(\pi\xb), & \kappa_{2}&= 1-2\sin(\pi\xb), \label{eq:1dDeltaProblemNL}\\
  \nu&= -\frac{3\pi}{\delta_{1}+\delta_{2}} &\mu&= 1 - 3\sin(\pi\xb) \;\;{\rm or} \;\; \mu=1. \nonumber
\end{align}
Note that both the proposed choices of $\mu$ are consistent with Theorem \ref{thm:delta-convergence}; however, while the first choice is exact (i.e. it corresponds to the difference of the manufactured local solution $u_1-u_2$), the second only matches the manufactured solutions' difference on $\Gamma_0$. For this reason, the latter is a more realistic choice, i.e. more likely to be used in practical settings. Our results will show that both choices guarantee the convergence to the local solution in the limit of vanishing nonlocality. 
We consider kernels of both type 1 and 2. Specifically, for $i=1,2$, we choose the fractional kernels $\gamma_i^F$ as in \eqref{eq:fractional-kernel} with $s_1=0.2$ and $s_2=0.4$ and the constant kernels $\gamma_i^C$ as in \eqref{eq:constant-kernel}.

\paragraph{$\delta$-Convergence illustrations in 1D}
As done in the previous section we consider three cases: all kernels are of type 1, all kernels are of type 2, and mixed kernels (i.e. constant kernels in $\Omega_1$ and fractional kernels in $\Omega_2$). Furthermore, we consider two different ratios between the horizons, specifically, \(\delta_{2}/\delta_{1}=2\) and \(\delta_{2}/\delta_{1}=1\) and study the convergence behavior as \(\delta_{1}\rightarrow 0\).

\begin{figure}
  \centering
  \includegraphics[width=0.6\textwidth]{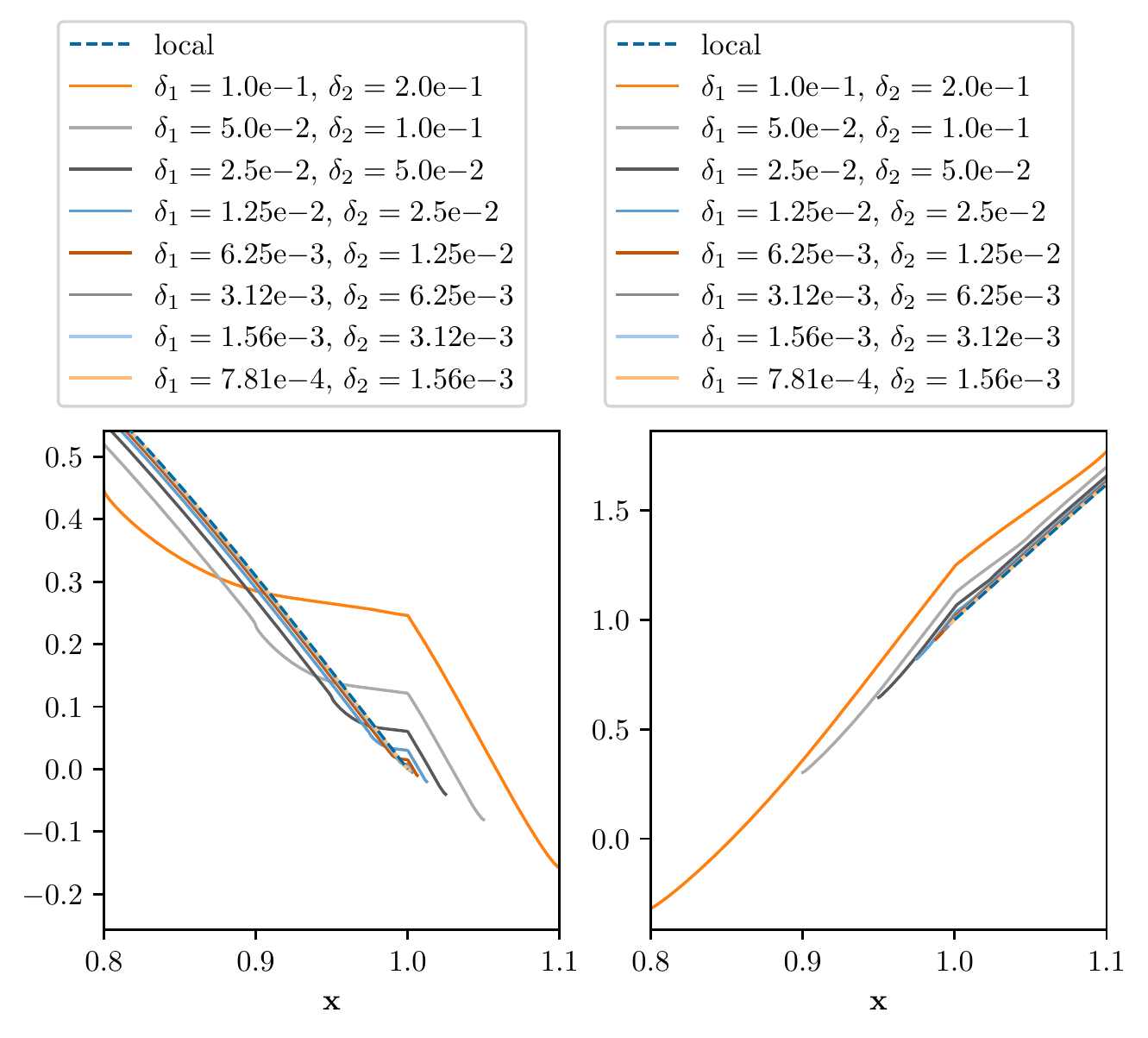} \\
  \includegraphics[width=0.6\textwidth]{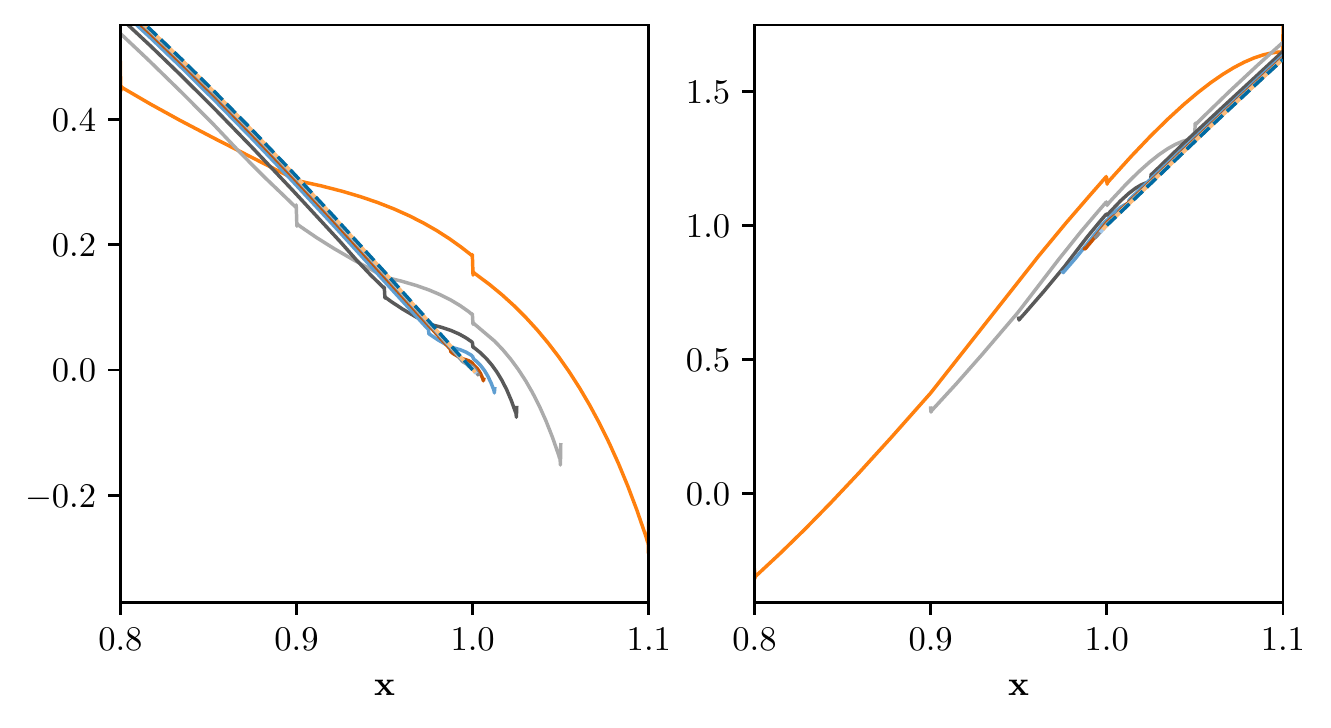}\\
  \includegraphics[width=0.6\textwidth]{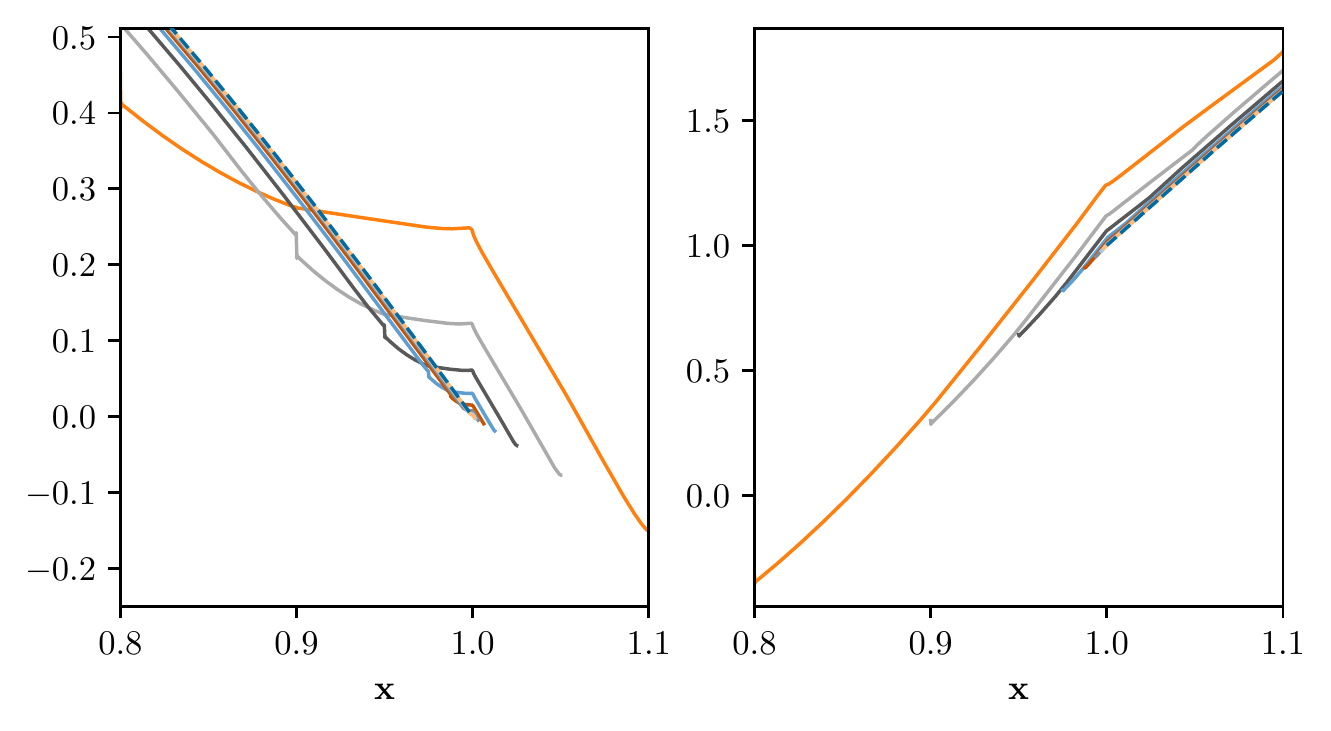}
  \caption{For $\mu=1-3\sin(\pi \xb)$, solutions $u_1$ (\emph{left}) and $u_2$ (\emph{right}) of the 1D problem \eqref{eq:1dDeltaProblemNL} near the interface, with \(\delta_{2}/\delta_{1}=2\) and \(h\approx 2\times10^{-4}\).
    \emph{Top:} fractional kernels on both subdomains with \(s_{1}=0.2\) and \(s_{2}=0.4\).
    \emph{Middle:} constant kernels on both subdomains.
    \emph{Bottom:} constant kernel on the left subdomain, fractional kernels on the right subdomain with \(s_{2}=0.4\).
  }
  \label{fig:deltaConvergenceVisual}
\end{figure}

\begin{figure}
  \centering
  \includegraphics[width=0.6\textwidth]{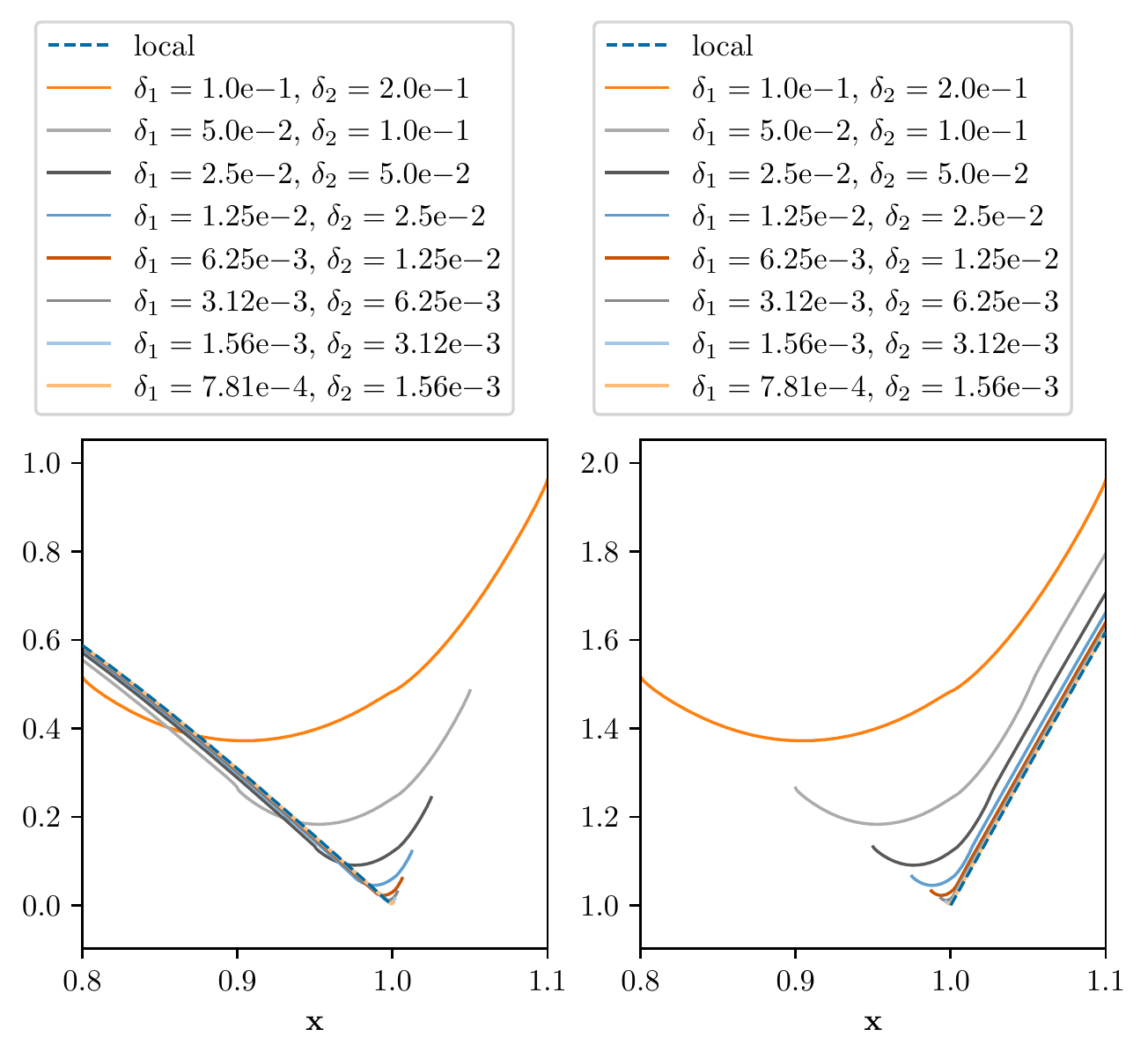} \\
    \includegraphics[width=0.6\textwidth]{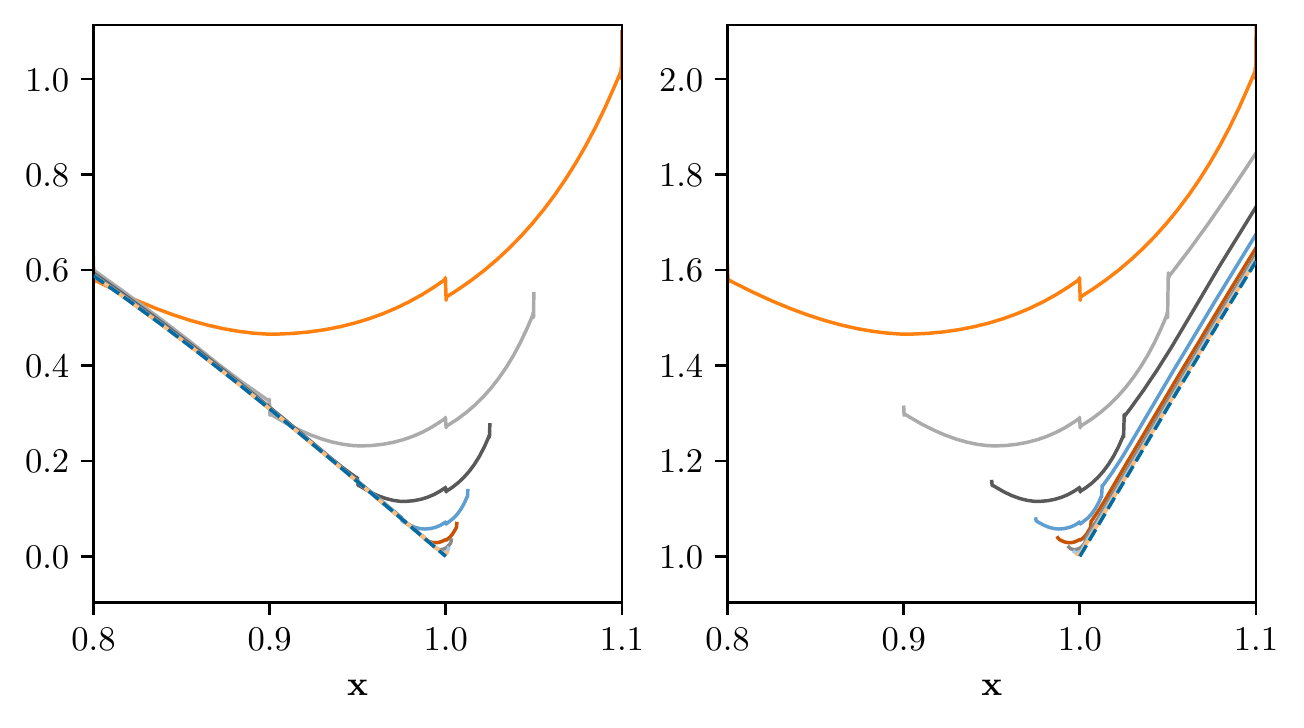}\\
  \includegraphics[width=0.6\textwidth]{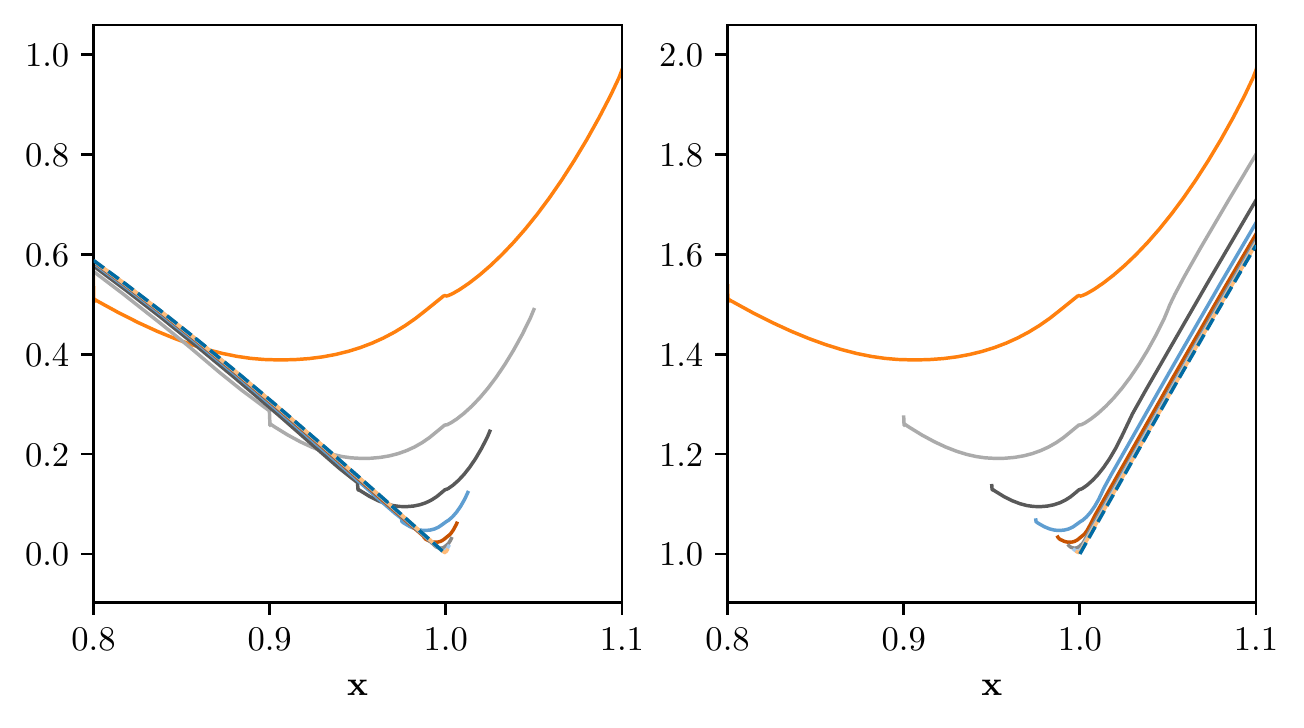}
  \caption{For $\mu=1$, solutions $u_1$ (\emph{left}) and $u_2$ (\emph{right}) of the 1D problem \eqref{eq:1dDeltaProblemNL} near the interface, with \(\delta_{2}/\delta_{1}=2\) and \(h\approx 2\times10^{-4}\).
    \emph{Top:} fractional kernels on both subdomains with \(s_{1}=0.2\) and \(s_{2}=0.4\).
    \emph{Middle:} constant kernels on both subdomains.
    \emph{Bottom:} constant kernel on the left subdomain, fractional kernels on the right subdomain with \(s_{2}=0.4\).}
  \label{fig:deltaConvergenceVisualmu1}
\end{figure}

We first provide a visual confirmation of the $\delta$-convergence for both choices of $\mu$ in Figures~\ref{fig:deltaConvergenceVisual} and \ref{fig:deltaConvergenceVisualmu1} respectively. Here, for \(\delta_{2}/\delta_{1}=2\) and all three combinations of kernel functions, we report the nonlocal solutions as $\delta_1\to 0$ and the exact local solution. We observe that at the limit of vanishing nonlocality the nonlocal solution matches the local one, as predicted by Theorem \ref{thm:delta-convergence}, regardless of the choice of $\mu$ in $\Gamma$, as long as $\mu|_{\Gamma_0}=m=1$. Thus, in all the remaining experiments, we use $\mu=1$ in $\Gamma$.

\begin{figure}
  \centering
  \includegraphics[width=0.32\textwidth]{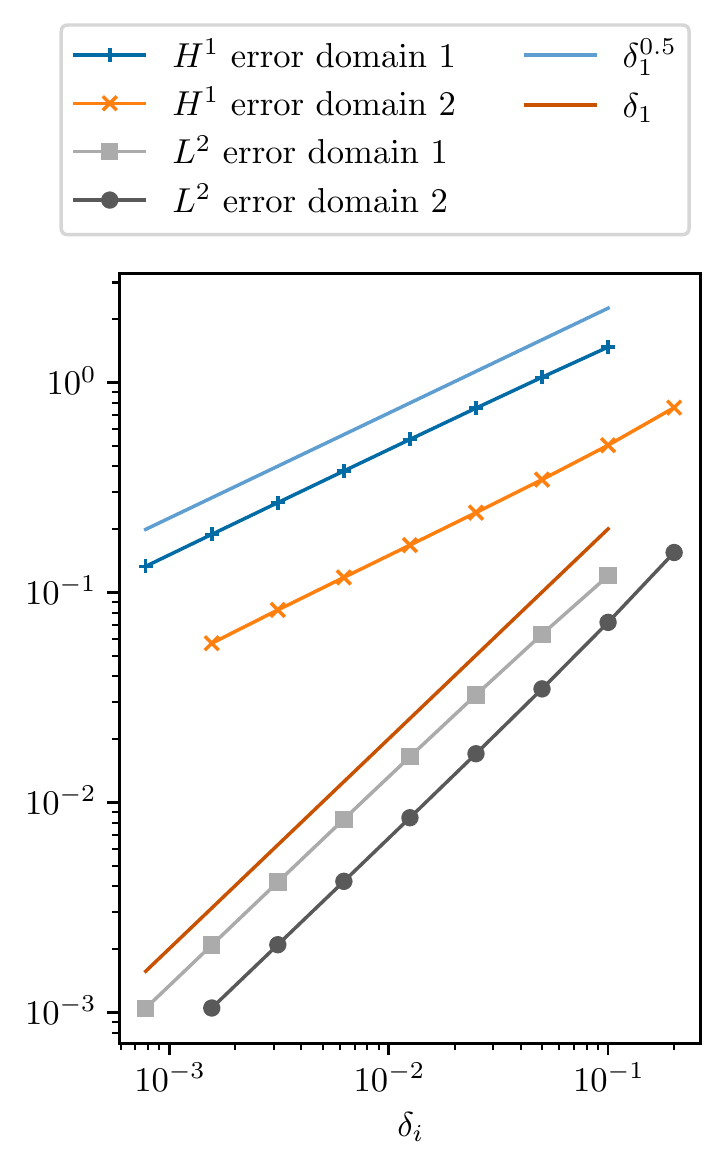}
  \includegraphics[width=0.32\textwidth]{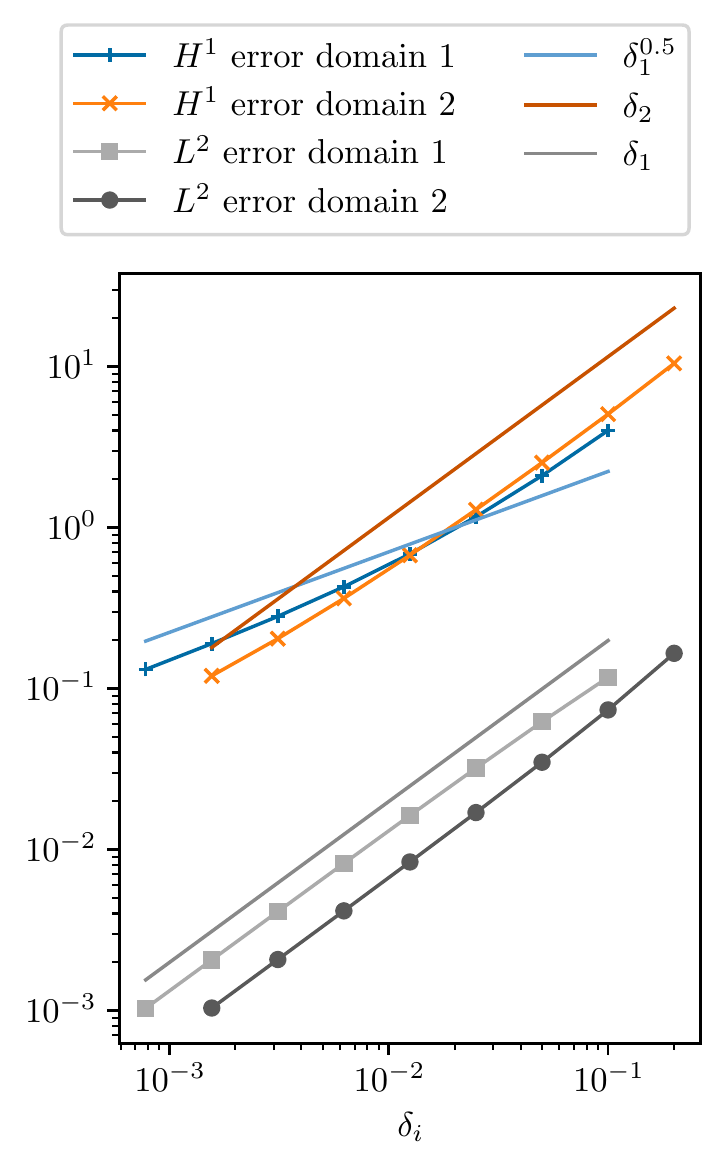}
  \includegraphics[width=0.32\textwidth]{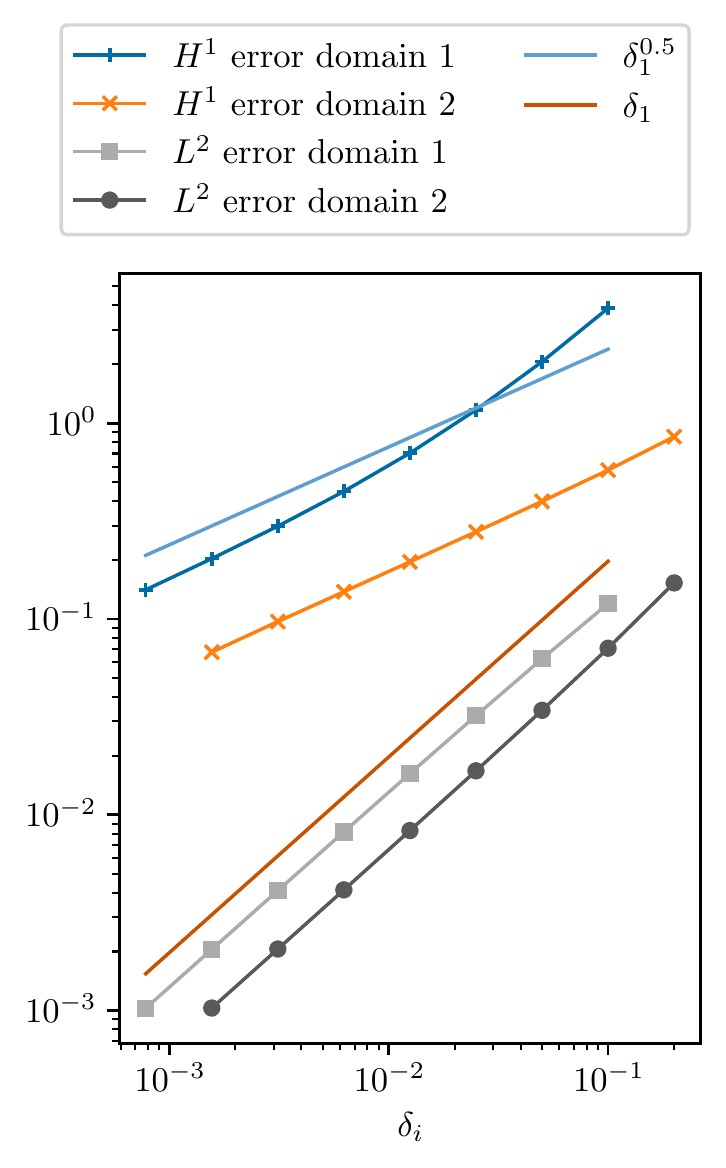}
  \caption{
    Convergence of the 1D problem \eqref{eq:1dDeltaProblemNL} with respect to the horizon, with \(\delta_{2}/\delta_{1}=2\) and \(h\approx 2\times10^{-4}\).
    \emph{Left:} fractional kernels on both subdomains with \(s_{1}=0.2\) and \(s_{2}=0.4\).
    \emph{Center:} constant kernels on both subdomains.
    \emph{Right:} constant kernel on the left subdomain, fractional kernels on the right subdomain with \(s_{2}=0.4\).
  }
  \label{fig:deltaConvergence2}
\end{figure}

\begin{figure}
  \centering
  \includegraphics[width=0.32\textwidth]{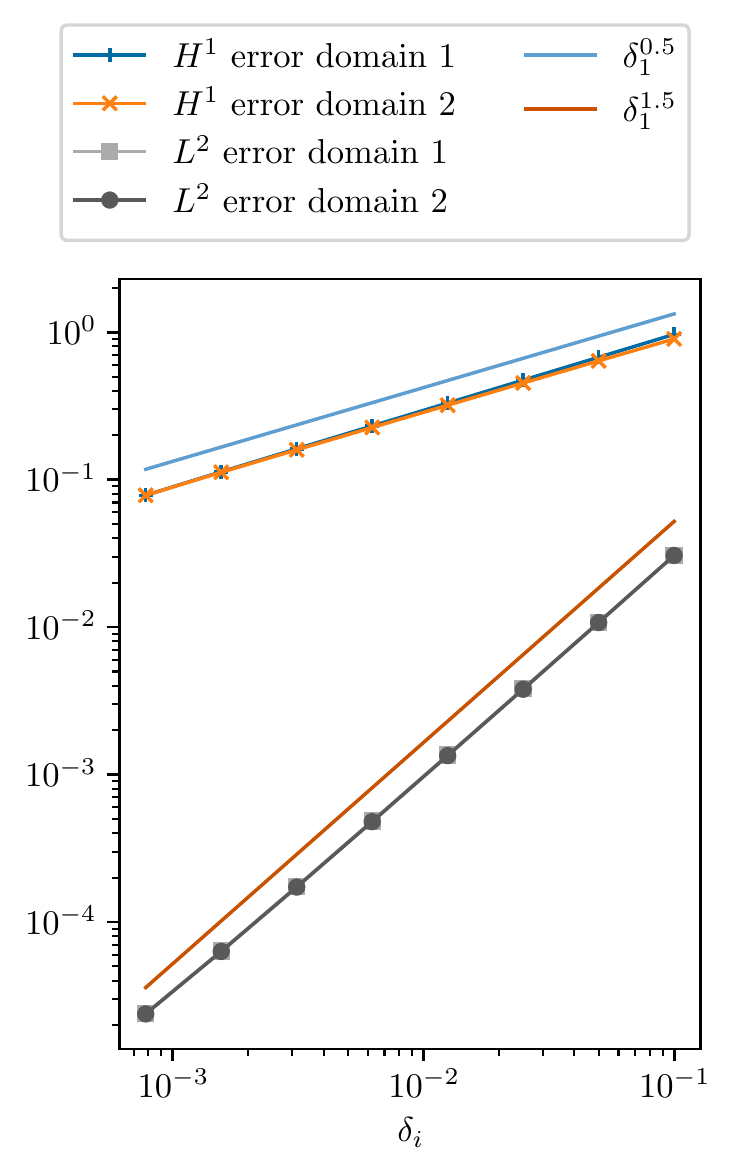}
  \includegraphics[width=0.32\textwidth]{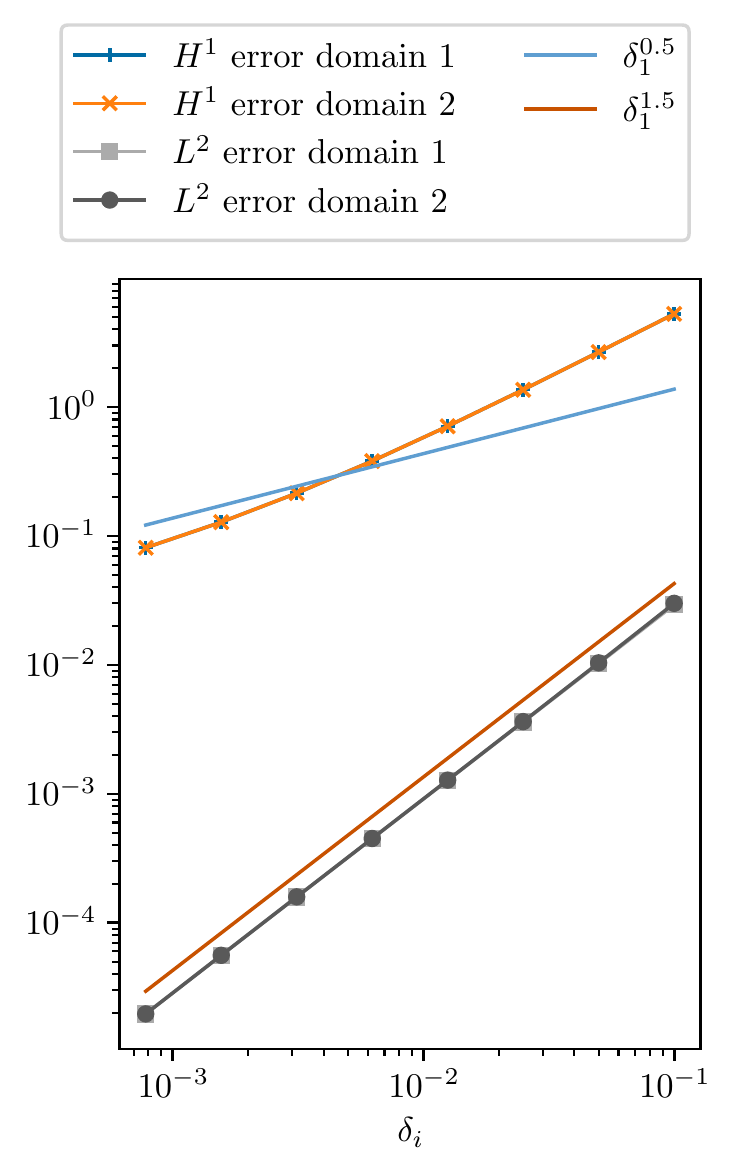}
  \includegraphics[width=0.32\textwidth]{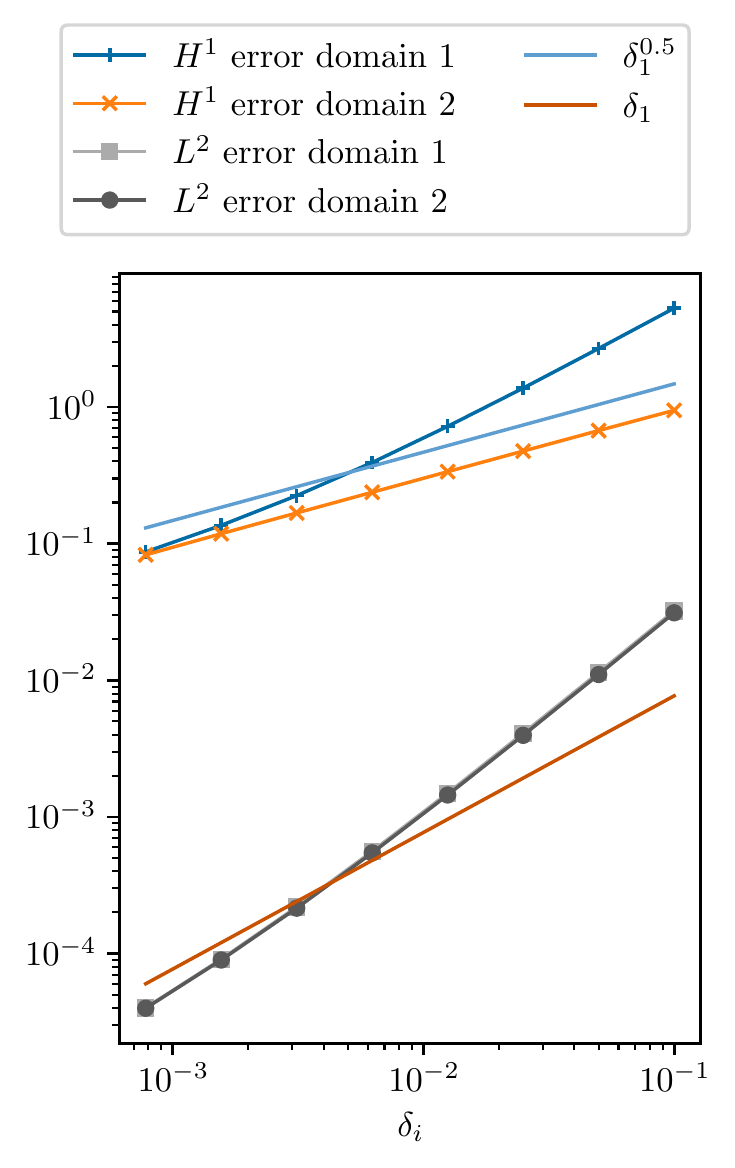}
  \caption{
    Convergence of the 1D problem \eqref{eq:1dDeltaProblemNL} with respect to the horizon, with \(\delta_{2}/\delta_{1}=1\) and \(h\approx 2\times10^{-4}\).
    \emph{Left:} fractional kernels on both subdomains with \(s_{1}=0.2\) and \(s_{2}=0.4\).
    \emph{Center:} constant kernels on both subdomains.
    \emph{Right:} constant kernel on the left subdomain, fractional kernels on the right subdomain with \(s_{2}=0.4\).
    }
  \label{fig:deltaConvergence1}
\end{figure}

Further confirmation of the $\delta$-convergence behavior and of the statement in Remark \ref{rem:conv-rate} is given by the results in Figures~\ref{fig:deltaConvergence2} and \ref{fig:deltaConvergence1}, where, for both horizon ratios (2 and 1 respectively) we study the $\delta$-convergence rate for all three combinations of kernels. In all runs, the mesh size is \(h \sim 2\times 10^{-4}\), this choice is small enough that the discretization error does not affect the $\delta$-convergence. Errors in the \(H^{1}(\Omega_{i})\)-seminorms and \(L^{2}(\Omega_{i})\)-norms are computed with respect to the known analytic solution of the local problem. We observe a $\delta$-convergence of order \(\mathcal{O}(\delta_{i}^{0.5})\) for the \(H^{1}\)-seminorm, as expected, and of order \(\mathcal{O}(\delta_{i})\) for the \(L^{2}\) norms. We also point out that for the case of ratio 1 (i.e. $\delta_1=\delta_2$) and kernels of identical types, the observed $\delta$-convergence rate in the $L^2(\Omega_i)$ norm is 1.5. This is possibly due to cancellations because the kernels are the same throughout the domain. These results also confirm the asymptotic compatibility of the proposed method.

\paragraph{Problem setting in 2D}
As in the previous two-dimensional experiments, we let $\Omega_{1}= (0,1)^{2}$ and $\Omega_{2}=(1,2)\times(0,1)$.
We first define the local interface problem and then choose the data of the nonlocal interface problem according to Theorem \ref{thm:delta-convergence}. Let the local source terms, the Dirichlet boundary conditions, and the jump terms be defined as
\begin{align*}
  f_{1}&= 10\pi^{2}\sin(\pi x_{1})\sin(2\pi x_{2}), & f_{2}&= 2\pi^{2}\sin(\pi x_{1})\sin(\pi x_{2})\\
  g_{1}&= 2+2\sin(\pi x_{1})\sin(2\pi x_{2}), & g_{2}&= 1-\sin(\pi x_{1})\sin(\pi x_{2})\\
  s&= -2\pi\sin(2\pi x_{2}) - \pi\sin(\pi x_{2}), &  m&=-1.
\end{align*}
The exact solution to the local problem is then given by $u_{1}= 2+2\sin(\pi x_{1})\sin(\pi x_{2})$ and $u_{2}= 1-\sin(\pi x_{1})\sin(\pi x_{2})$. The corresponding nonlocal source terms, Dirichlet volume constraints, and jump terms are then given by
\begin{align}
  \zeta_{1}&= 10\pi^{2}\sin(\pi x_{1})\sin(2\pi x_{2}), & \zeta_{2}&= 2\pi^{2}\sin(\pi x_{1})\sin(\pi x_{2}) \nonumber\\
  \kappa_{1}&= 2+2\sin(\pi x_{1})\sin(2\pi x_{2}), & \kappa_{2}&= 1-\sin(\pi x_{1})\sin(\pi x_{2}),\label{eq:2dDeltaProblemNL}\\
  \nu&= -\frac{2\pi\sin(2\pi x_{2}) + \pi\sin(\pi x_{2})}{\delta_{1}+\delta_{2}} &\mu&= -1. \nonumber
\end{align}
We consider type 1 and type 2 kernels with normalization as given in \eqref{eq:fractional-kernel2d} and \eqref{eq:constant-kernel2d}.

\begin{figure}
  \centering
  \includegraphics[width=0.32\textwidth]{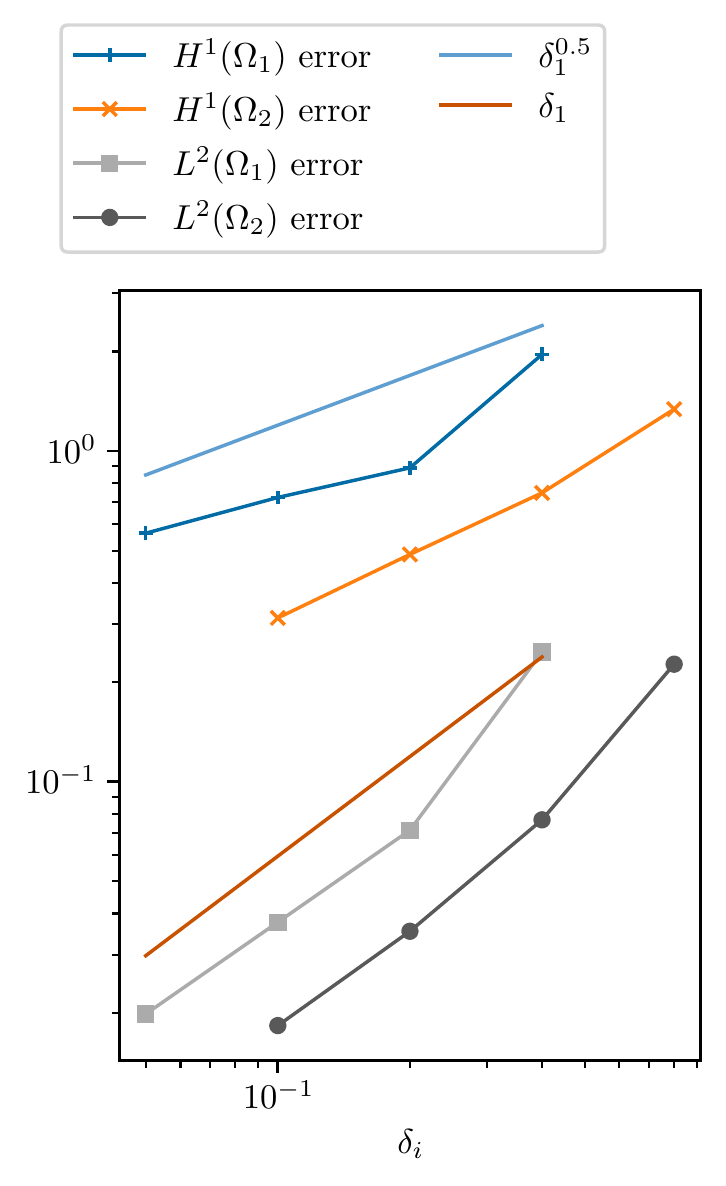}
  \includegraphics[width=0.32\textwidth]{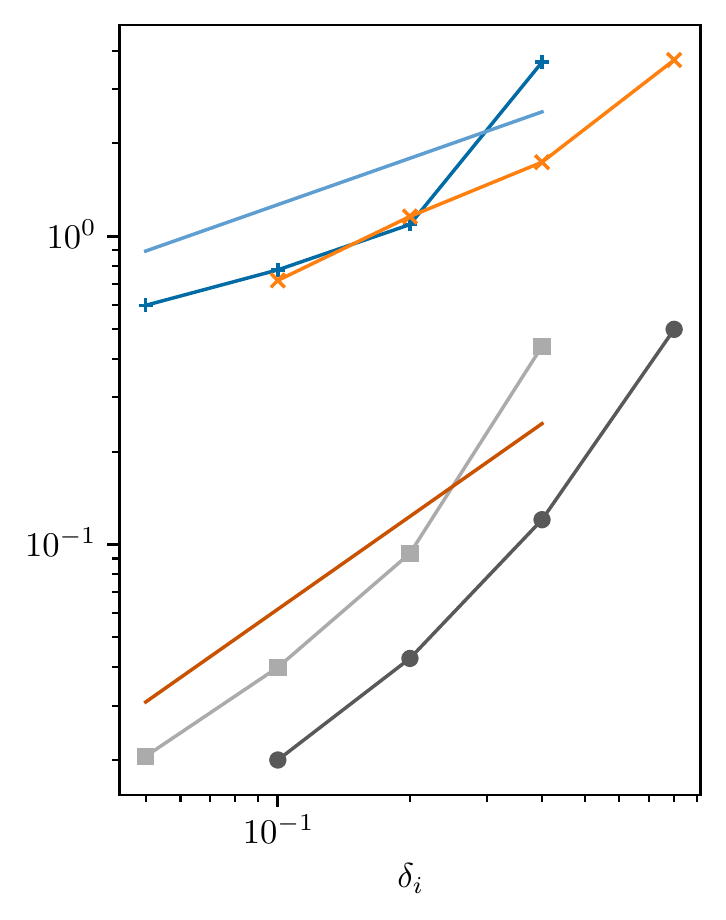}
  \includegraphics[width=0.32\textwidth]{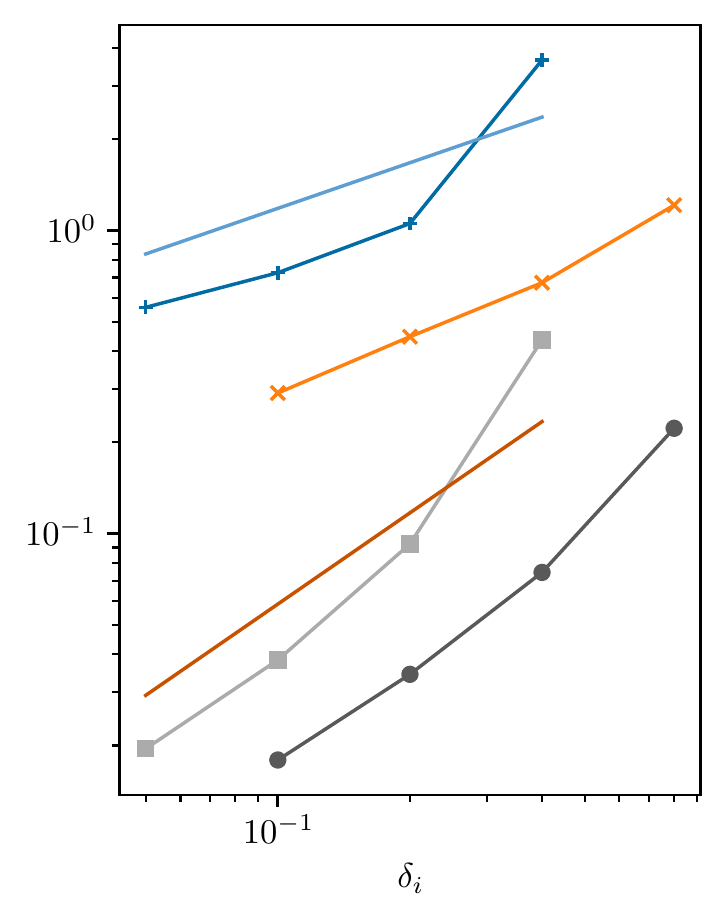}
  \caption{
    Convergence of the 2D problem \eqref{eq:2dDeltaProblemNL} with respect to the horizon, with \(\delta_{2}/\delta_{1}=2\) and \(h\approx \delta_{1}/4\).
    \emph{Left:} fractional kernels on both subdomains with \(s_{1}=0.2\) and \(s_{2}=0.4\).
    \emph{Center:} constant kernels on both subdomains.
    \emph{Right:} constant kernel on the left subdomain, fractional kernels on the right subdomain with \(s_{2}=0.4\).
  }
  \label{fig:deltaConvergence2D2}
\end{figure}

\begin{figure}
  \centering
  \includegraphics[width=0.32\textwidth]{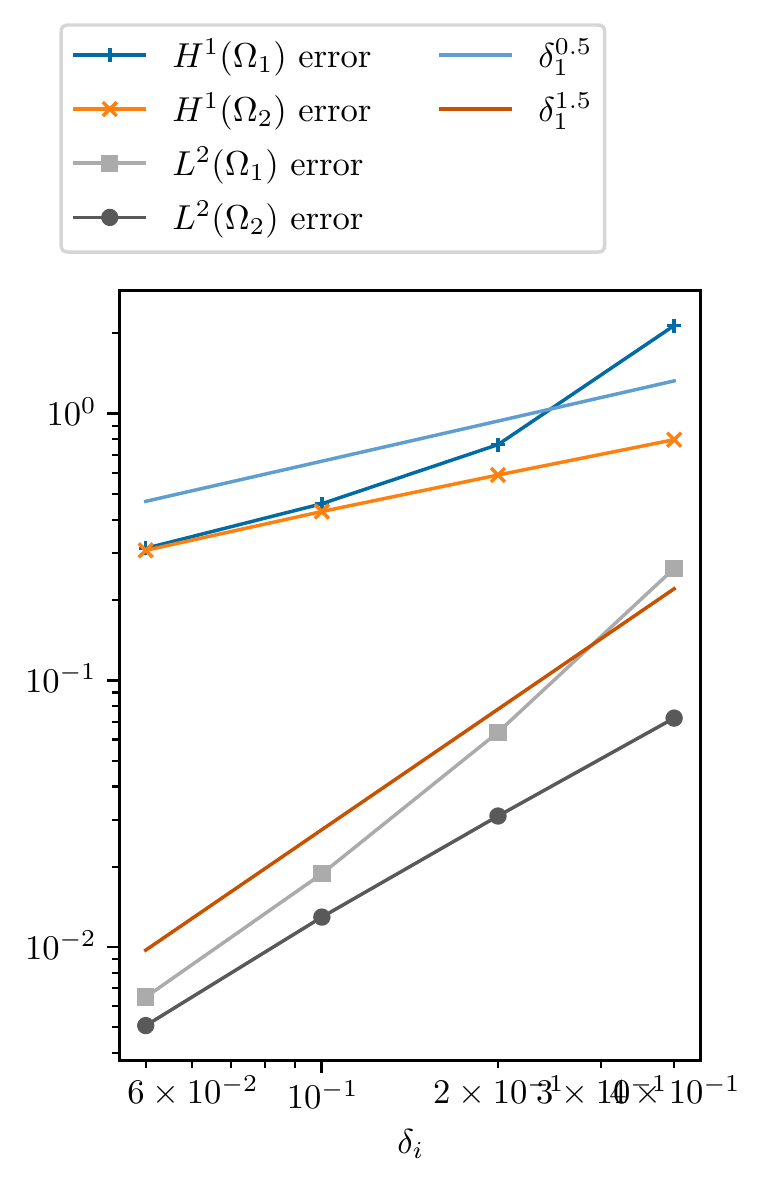}
  \includegraphics[width=0.32\textwidth]{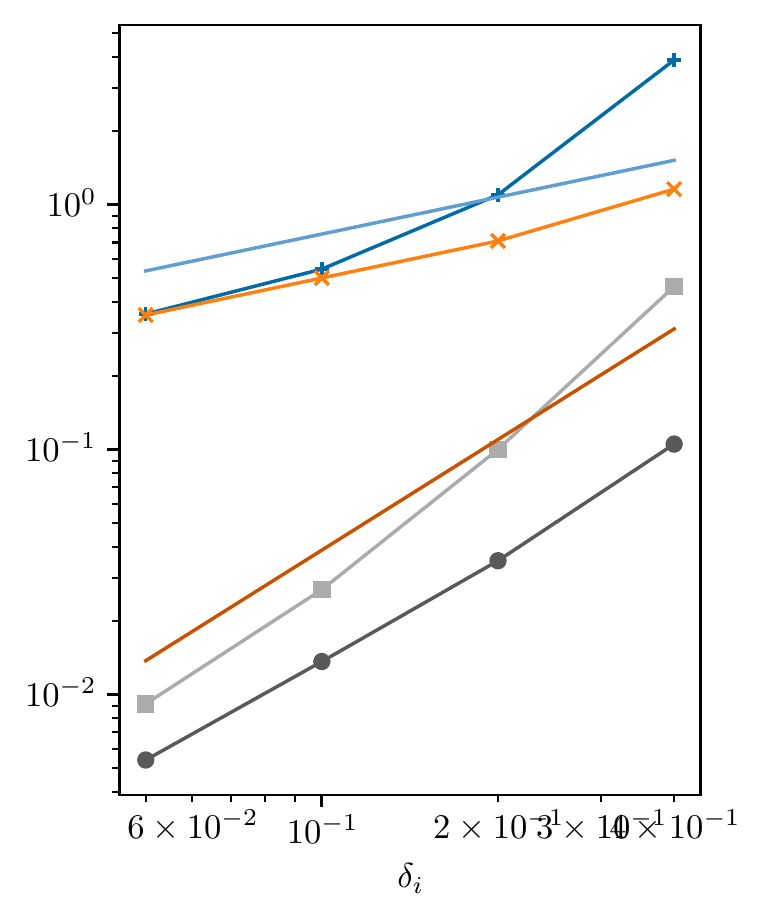}
  \includegraphics[width=0.32\textwidth]{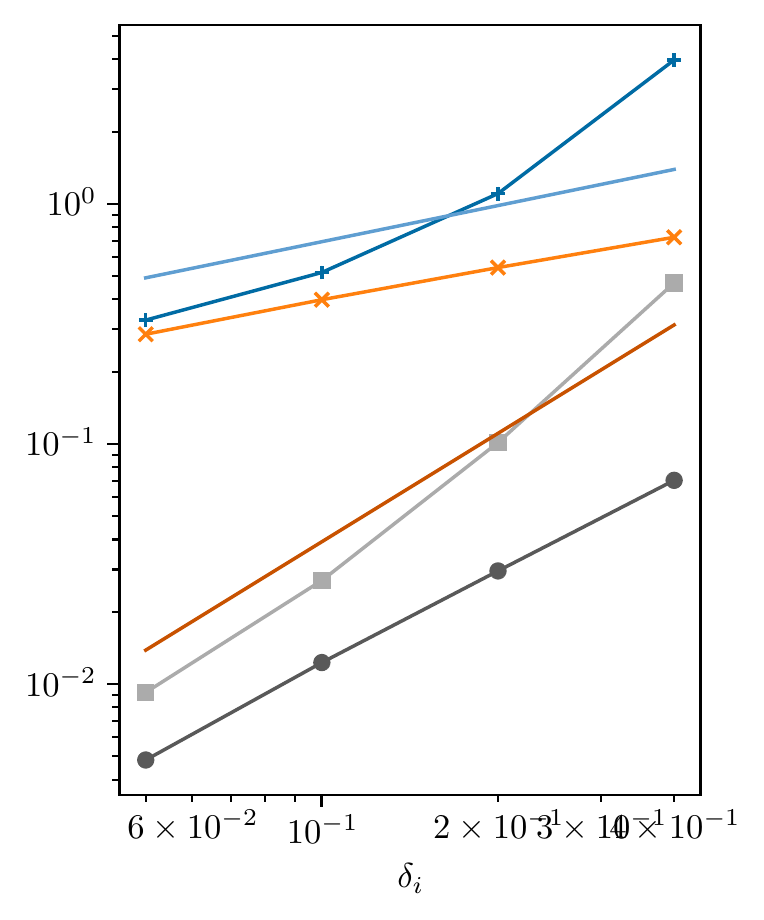}
  \caption{
    Convergence of the 2D problem \eqref{eq:2dDeltaProblemNL} with respect to the horizon, with \(\delta_{2}/\delta_{1}=1\) and \(h\approx \delta_{1}/4\).
    \emph{Left:} fractional kernels on both subdomains with \(s_{1}=0.2\) and \(s_{2}=0.4\).
    \emph{Center:} constant kernels on both subdomains.
    \emph{Right:} constant kernel on the left subdomain, fractional kernels on the right subdomain with \(s_{2}=0.4\).
    }
  \label{fig:deltaConvergence2D1}
\end{figure}

\paragraph{$\delta$-Convergence illustrations in 2D}
As with the one-dimensional experiments, we consider two ratios of horizon, \(\delta_{2}/\delta_{1}=2\) and \(\delta_{2}/\delta_{1}=1\) and explore the convergence of the nonlocal solutions to the local ones as \(\delta_{1}\rightarrow0\) in Figures~\ref{fig:deltaConvergence2D2} and \ref{fig:deltaConvergence2D1}. We select the mesh size as \(h\approx \delta_{1}/4\), and compute errors in \(H^{1}(\Omega_{i})\)-seminorm and \(L^{2}(\Omega_{i})\)-norm with respect to the known solution of the local problem.
Also in this case, we observe convergence of order \(\mathcal{O}(\delta_{i}^{0.5})\) in \(H^{1}\)-seminorm and at least of order \(\mathcal{O}(\delta_{i})\) in \(L^{2}\)-norm.
Again, for a horizon ratio of \(\delta_{2}/\delta_{1}=1\), we observe an improved rate of \(\mathcal{O}(\delta_{i}^{1.5})\) in \(L^{2}\)-norm.
In contrast to the one-dimensional experiments, this behavior can also be observed in the case of type 2 kernel in $\Omega_1$ and type 1 kernel in $\Omega_2$.

\section{Concluding remarks}\label{sec:conclusion}
In this work we introduced, for the first time, a mathematically rigorous formulation for nonlocal interface problems that feature solution and/or flux jumps at the interface. The proposed weak form of the nonlocal interface problem is well-posed and mimics the weak formulation of the corresponding local interface problem with jumps. Furthermore, when the nonlocal data are appropriately prescribed, our analysis shows that the nonlocal solution approaches the solution of the local problem (see Theorem \ref{thm:delta-convergence}). For the proposed weak form, we also introduce a finite element discretization and show that, for piecewise linear finite element spaces, the nonlocal solution features an optimal, quadratic convergence rate in the $L^2$ metric. Furthermore, we show that the numerical solution converges to an analytical solution of the corresponding local problem as the horizon vanishes, confirming the asymptotic compatibility of the proposed scheme. 
Finally, we point out that our reinterpretation of the proposed formulation as a single-domain problem with a piecewise defined kernel makes the implementation task trivial. In fact, the presence of discontinuities in the model parameters only requires a re-definition of the single-domain kernel, and the presence of the jumps simply corresponds to additional forcing terms applied on and around the nonlocal interface. In other words, the structure of a nonlocal finite element code is the same as the one of a standard, single-domain problem. 

Naturally, extensions of this work include demonstrating the validity and efficacy of the nonlocal interface formulation on relevant mechanics problems. Thus, we plan to extend the current formulation to, e.g., the state-based peridynamic model \cite{silling2007peridynamic}. Moreover, though not considered in this paper, the treatment of multiple interfaces can be conducted following the same formulation and the same choice of interface data. Thus, the implementation of this more realistic setting is straightforward.

\section*{Acknowledgments}
M. D'Elia and P. Bochev were partially supported by the U.S. Department of Energy, Office of Advanced Scientific Computing Research under the Collaboratory on Mathematics and Physics-Informed Learning Machines for Multiscale and Multiphysics Problems (PhILMs) project. P. Bochev, M. D'Elia, C. Glusa, and M. Gunzburger were also supported by the Sandia National Laboratories Laboratory-directed Research and Development (LDRD) program, project 218318. Sandia National Laboratories is a multimission laboratory managed and operated by National Technology and Engineering Solutions of Sandia, LLC., a wholly owned subsidiary of Honeywell International, Inc., for the U.S. Department of Energy's National Nuclear Security Administration contract number DE-NA0003525. This paper, SAND2022-2944 O, describes objective technical results and analysis. Any subjective views or opinions that might be expressed in the paper do not necessarily represent the views of the U.S. Department of Energy or the United States Government.
\bibliographystyle{abbrv}
\bibliography{references}

\end{document}